\documentclass[12pt]{article}
\usepackage{amsfonts}
\usepackage{bbm}
\usepackage{mathrsfs}
\usepackage{amssymb}

\def\hang{\hangindent\parindent}
\def\textindent#1{\indent\llap{#1\enspace}\ignorespaces}
\def\re{\par\hang\textindent}

\title{Computation of  Minimal Graded Free
Resolutions over\\ $\mathbb{N}$-Graded Solvable Polynomial Algebras}

\author{Huishi Li\\
{\small Department of Applied Mathematics},
{\small College of Information Science and Technology}\\
{\small Hainan University}, {\small  Haikou 570228, China}}

\date{}

\begin{document}
\maketitle
\begin{center}
\begin{minipage}{140mm}
{\small {\bf Abstract.} It is shown that the methods and algorithms,
developed in (A. Capani et al., Computing minimal finite free
resolutions, {\it Journal of Pure and Applied Algebra}, (117\&
118)(1997), 105 -- 117; M. Kreuzer and L. Robbiano, {\it
Computational Commutative Algebra 2}, Springer, 2005.) for computing
minimal homogeneous generating sets  of graded submodules and graded
quotient modules of  free modules  over a commutative polynomial
algebra, can be adapted for computing minimal homogeneous generating
sets  of graded submodules and graded quotient modules of free
modules over a weighted $\mathbb{N}$-graded solvable polynomial
algebra, where solvable polynomial algebras are in the sense of
(A.~Kandri-Rody and V.~Weispfenning, Non-commutative Gr\"obner bases
in algebras of solvable type. {\it J. Symbolic Comput.}, 9(1990),
1--26). Consequently, algorithmic procedures for computing minimal
finite graded free resolutions over weighted $\mathbb{N}$-graded
solvable polynomial algebras are achieved. }
\end{minipage}\end{center} {\parindent=0pt\vskip 6pt

{\bf 2010 Mathematics subject classification} Primary 16W70; 
Secondary 16Z05.\vskip 6pt

{\bf Key words} Solvable polynomial algebra,  homogeneous generating 
set, Gr\"obner basis, free resolution.}

\def\v5{\vskip .5truecm}\def\QED{\hfill{$\Box$}}\def\hang{\hangindent\parindent}
\def\textindent#1{\indent\llap{#1\enspace}\ignorespaces}
\def\item{\par\hang\textindent}
\def \r{\rightarrow}\def\OV#1{\overline {#1}}
\def\normalbaselines{\baselineskip 24pt\lineskip 4pt\lineskiplimit 4pt}
\def\mapdown#1{\llap{$\vcenter {\hbox {$\scriptstyle #1$}}$}
                                \Bigg\downarrow}
\def\mapdownr#1{\Bigg\downarrow\rlap{$\vcenter{\hbox
                                    {$\scriptstyle #1$}}$}}
\def\mapright#1#2{\smash{\mathop{\longrightarrow}\limits^{#1}_{#2}}}
\def\NZ{\mathbb{N}}

\def\LH{{\bf LH}}\def\LM{{\bf LM}}\def\LT{{\bf
LT}}\def\KX{K\langle X\rangle} \def\KS{K\langle X\rangle}
\def\B{\mathscr{B}} \def\LC{{\bf LC}} \def\G{{\cal G}} \def\FRAC#1#2{\displaystyle{\frac{#1}{#2}}}
\def\SUM^#1_#2{\displaystyle{\sum^{#1}_{#2}}} \def\O{{\cal O}}  \def\J{{\bf J}}\def\BE{\B (e)}
\def\PRCVE{\prec_{\varepsilon\hbox{-}gr}}\def\BV{\B (\varepsilon )}\def\PRCEGR{\prec_{\mathbbm{e}\hbox{\rm -}gr}}

\def\KS{K\langle X\rangle}
\def\LR{\langle X\rangle}
\def\HL{{\rm LH}}\def\NB{\mathbb{N}}

\vskip 1truecm

\section*{1. Introduction}
Since the late 1980s, the  Gr\"obner basis theory for commutative 
polynomial algebras and their modules (cf. [Bu1, 2],  [Sch], [BW], 
[AL2], [Fr\"ob], [KR1, 2]) has been successfully generalized to 
(noncommutative) solvable polynomial algebras and their modules (cf. 
[AL1], [Gal], [K-RW], [Kr2], [LW], [Li1], [Lev]). It is now well 
known that the class of solvable polynomial algebras covers numerous 
significant algebras such as enveloping algebras of Lie algebras, 
Weyl algebras (including  algebras of partial differential operators 
with polynomial coefficients over a field $K$ of characteristic 0), 
more generally a large number of operator algebras, iterated Ore 
extensions, and many quantum (quantized) algebras. In particular, 
after [K-RW] successfully established a noncommutative version of 
the Buchberger's criterion and  a noncommutative version of 
Buchberger algorithm for computing (one-sided, two-sided) Gr\"obner 
bases of (one-sided, two-sided) ideals in general solvable 
polynomial algebras (see the module versions presented as Theorem 
2.5 and {\bf Algorithm 1} in the current paper), the noncommutative 
version of Buchberger algorithm  has been implemented in some 
well-developed computer algebra systems, such as \textsc{Modula-2} 
[KP] and \textsc{Singular} [DGPS]. Based on such an {\it effective 
Gr\"obner basis theory} and the fact that {\it every solvable 
polynomial algebra is a} ({\it left and right}) {\it Noetherian 
domain of finite global homological dimension} (see Theorem 2.3  and 
Theorem 5.3  in the subsequent sections), in this paper we show that 
the    methods and algorithms, developed in  ([CDNR], [KR]) for 
computing minimal homogeneous generating sets  of graded submodules 
and graded quotient modules of free modules over a commutative 
polynomial algebra, can be adapted for  computing minimal 
homogeneous generating sets  of graded submodules and  graded 
quotient modules of free modules over a weighted $\mathbb{N}$-graded 
(noncommutative) solvable polynomial algebra $A=K[a_1,\ldots ,a_n]$ 
with the degree-0 homogeneous part $A_0=K$, where $K$ is a field, 
and consequently, algorithmic procedures for computing minimal 
finite graded free resolutions over weighted $\mathbb{N}$-graded 
solvable polynomial algebras can be achieved.  More precisely, after 
the preliminary Section 2, we present algorithms in Section 3 for 
computing $n$-truncated left Gr\"obner bases and minimal homogeneous 
generating sets of graded submodules in free modules over  $A$; in 
Section 4 we present an algorithm for computing a minimal 
homogeneous generating set of a graded quotient module $M=L/N$  of a 
free  $A$-module $L$; and in the final Section 5 we present 
algorithmic  procedures for computing a minimal finite graded free 
resolution of a finitely generated graded $A$-module $M\cong L/N$.
\par

Noticing that commutative polynomial algebras are certainly the type 
of $\NZ$-graded solvable polynomial algebras we specified, and that 
the noncommutative version of Buchberger's criterion as well as 
Buchberger algorithm for modules over solvable polynomial algebras 
(Theorem 2.5 and {\bf Algorithm 1} presented in the end of Section 
1) looks as if working the same way as in the commutative case by 
reducing the S-polynomials, one might think that the extension of 
methods and algorithms provided by ([CDNR], [KR2]) to  modules over 
noncommutative $\mathbb{N}$-graded solvable polynomial algebras 
could be naturally holding true as a folklore. However, from the 
literature  (e.g. [AL2], [BW], [Eis], [Fr\"ob], [KR1,2]) we learnt 
that in developing the Gr\"obner basis (including the $n$-truncated 
Gr\"obner basis) theory for a commutative polynomial $K$-algebra 
$R=K[x_1,\ldots ,x_n]$,  two featured algebraic structures play the 
key role in both the theoretical proofs and technical calculations, 
namely {\parindent=.5truecm\par

\item{$\bullet$} the  multiplicative
monoid $\B$, where $\B$ is the PBW basis of $R$, which is 
furthermore turned into an ordered multiplicative monoid with 
respect to any monomial ordering $\prec$; \par

\item{$\bullet$} monomial ideals, i.e., ideals
generated by monomials from $\B$.\par}{\parindent=0pt\par

For instance,  a version of Dickson's lemma for monomial ideals 
holds true, thereby a Gr\"obner basis of an ideal $I$ in $R$ is 
usually defined (or characterized) in terms of generators of the 
monomial ideal generated by leading monomials of $I$; in the proof 
of Buchberger's criterion, reduction of a monomial does not cause 
any trouble (e.g. see [AL2], P.41, l -5); especially,  the already 
known Noetherianess of $R$ (or Dickson's lemma for monomial ideals) 
guarantees the termination of Buchberger algorithm (e.g. see [AL2], 
P.43), and this algorithm,  in turn, gives rise to more relevant 
algorithms not only for ideals but also for submodules of free 
$R$-modules (e.g. [KR2], Proposition 4.5.10, Theorem 4.6.3). While 
due to the {\it noncommutativity} of a solvable polynomial algebra 
$A=K[a_1,\ldots ,a_n]$, {\it the PBW basis $\B$ of $A$ is  no longer 
a multiplicative monoid}. Thereby, in developing a (one-sided, 
two-sided) Gr\"obner basis theory of $A$, all jobs using reduction 
by monomials from $\B$ cannot be simply replicated from the 
commutative case, (one-sided, two-sided) ideals generated by 
monomials from $\B$ can no longer  play the  role  as in the 
commutative case, and the Noetherianess of $A$ is not  known until 
the existence of finite Gr\"obner bases for (one-sided) ideals is 
algorithmically established (note that in general $A$ is not 
necessarily an iterated Ore extension of the base ring $K$ or some 
Noetherian ring). Since  {\bf Algorithm 2} and {\bf Algorithm 3}  to 
be presented in Section 3 essentially depend on {\bf Algorithm 1} 
presented in the next section, at this point, one is referred to 
[K-RW] for the  nontrivial and detailed argumentation on how the 
barrier of noncommutativity is broken down, in order to reach the 
main results as we recalled in Section 2 (Theorem 2.3, Theorem 2.4, 
Theorem 2.5, {\bf Algorithm 1} though this is for modules). 
Moreover, so far in the literature there had been no a clear and 
systematical presentation  showing that the  commutative 
$n$-truncated Gr\"obner basis theory and the algorithmic principle 
for $\NZ$-graded modules presented in ([CDNR], [KR2]) can be adapted 
for $\NZ$-graded modules over general noncommutative $\NZ$-graded 
solvable polynomial $K$-algebras with the degree-0 homogeneous part 
equal to $K$. So, from a mathematical viewpoint, we are naturally 
concerned about how to trust that this is a true story and then, how 
to give a precise quotation source when the relevant results are  
applied to other noncommutative cases (for instance, in [Li5]). 
Following the rule of ``{\it to see is to believe}", which we 
understand as understanding more than merely observing, all what we 
pointed out above motivates us to provide a detailed argumentation 
and demonstration on the topic of this paper, to which one may also 
compare with the corresponding argumentations given in ([KR2], 
Chapter 4).}\par

In the literature, a finitely generated  $\NZ$-graded $K$-algebra 
$A=\oplus_{p\in\NZ}A_p$ with the degree-0 homogeneous part $A_0=K$ 
is referred to as a {\it connected $\NZ$-graded $K$-algebra}. 
Concerning introductions to minimal resolutions of graded modules 
over a (commutative or noncommutative) connected $\NZ$-graded 
$K$-algebra (or more generally an $\NZ$-graded local $K$-algebra) 
and relevant results, one may  refer to ([Eis], Chapter 19), ([Kr1], 
Chapter 3), and [Li3].\par

Throughout  this paper, $K$ denotes a field, $K^*=K-\{ 0\}$;
$\mathbb{N}$ denotes the additive monoid of all nonnegative
integers, and $\mathbb{Z}$ denotes the additive group of all
integers; all algebras are associative $K$-algebras with the
multiplicative identity 1, and modules over an algebra are meant
left unitary modules. \v5

\section*{2. Preliminaries}
In this section we  recall briefly some basics on Gr\"obner basis
theory for solvable polynomial algebras and their modules. The main
references are [AL1], [Gal], [K-RW], [Kr], [LW], [Li1], [Li4] and 
[Lev]. \v5

Let $K$ be a field and let $A=K[a_1,\ldots ,a_n]$ be a finitely
generated $K$-algebra with the {\it minimal set of generators} $\{
a_1,\ldots ,a_n\}$. If, for some permutation $\tau =i_1i_2\cdots
i_n$ of $1,2,\ldots ,n$, the set $\B =\{
a^{\alpha}=a_{i_1}^{\alpha_1}\cdots a_{i_n}^{\alpha_n}~|~\alpha
=(\alpha_1,\ldots ,\alpha_n)\in\NZ^n\} ,$ forms a $K$-basis of $A$,
then $\B$ is referred to as a {\it PBW $K$-basis} of $A$. It is
clear that if $A$ has a PBW $K$-basis, then we can always assume
that $i_1=1,\ldots ,i_n=n$. Thus, we make the following convention
once for all.{\parindent=0pt\v5

{\bf Convention} From now on in this paper, if we say that an
algebra $A$ has the PBW $K$-basis $\B$, then it means that
$$\B =\{ a^{\alpha}=a_{1}^{\alpha_1}\cdots
a_{n}^{\alpha_n}~|~\alpha =(\alpha_1,\ldots ,\alpha_n)\in\NZ^n\} .$$
Moreover, adopting the commonly used terminology in computational
algebra, elements of $\B$ are referred to as {\it monomials} of
$A$.} \v5

Suppose that $A$ has the PBW $K$-basis $\B$ as presented above and
that  $\prec$ is a total ordering on $\B$. Then every nonzero
element $f\in A$ has a unique expression
$$f=\lambda_1a^{\alpha (1)}+\lambda_2a^{\alpha (2)}+\cdots +\lambda_ma^{\alpha (m)},~
\lambda_j\in K^*,~a^{\alpha
(j)}=a_1^{\alpha_{1j}}a_2^{\alpha_{2j}}\cdots a_n^{\alpha_{nj}}\in\B
,~1\le j\le m.$$ If  $a^{\alpha (1)}\prec a^{\alpha (2)}\prec\cdots
\prec a^{\alpha (m)}$ in the above representation, then the {\it
leading monomial of $f$} is defined as $\LM (f)=a^{\alpha (m)}$, the
{\it leading coefficient of $F$} is defined as $\LC (f)=\lambda_m$,
and the {\it leading term of $f$} is defined as $\LT
(f)=\lambda_ma^{\alpha (m)}$. {\parindent=0pt\v5

{\bf 2.1. Definition}  Suppose that the $K$-algebra $A=K[a_1,\ldots
,a_n]$ has the PBW $K$-basis $\B$. If $\prec$ is a total ordering on
$\B$ that satisfies the following three
conditions:}{\parindent=1.34truecm\par

\item{(1)} $\prec$ is a well-ordering;\par

\item{(2)} For $a^{\gamma},a^{\alpha},a^{\beta}, a^{\eta}\in\B$, if
$a^{\alpha}\prec a^{\beta}$ and $\LM
(a^{\gamma}a^{\alpha}a^{\eta})$, $\LM
(a^{\gamma}a^{\beta}a^{\eta})\not\in K$, then $\LM
(a^{\gamma}a^{\alpha}a^{\eta})\prec\LM
(a^{\gamma}a^{\beta}a^{\eta})$;\par

\item{(3)} For $a^{\gamma},a^{\alpha},a^{\beta},a^{\eta}\in\B$, if $a^{\beta}\ne
a^{\gamma}$, and $a^{\gamma}=\LM (a^{\alpha}a^{\beta}a^{\eta})$,
then $a^{\beta}\prec a^{\gamma}$ (thereby $1\prec a^{\gamma}$ for
all $a^{\gamma}\ne 1$),\par}{\parindent=0pt

then $\prec$ is called a {\it monomial ordering} on $\B$ (or a
monomial ordering on $A$). }\par

If $\prec$ is a monomial ordering on $\B$, then we call $(\B
,\prec)$ an {\it admissible system} of $A$.{\parindent=0pt\v5

{\bf Remark.} (i) Definition 2.1 is indeed borrowed from the theory 
of Gr\"obner bases  for  general finitely generated $K$-algebras, in 
which the algebras considered may  be noncommutative, may have 
divisors of zero, and the $K$-bases used may  not be a PBW basis, 
but with a (one-sided, two-sided) monomial ordering such algebras 
may theoretically have a (one-sided, two-sided) Gr\"obner basis 
theory. For more details on this topic, one may referrer to ([Li2], 
Section 3.1 of Chapter 3 and Section 8.3 of Chapter 8). Also, to see 
the essential difference between Definition 2.1 and the classical 
definition of a monomial ordering in the commutative case, one may 
refer to (Definition 1.4.1 and the proof of Theorem 1.4.6 given in 
[AL2]). \vskip 6pt

(ii)  Note that the conditions (2) and (3) in Definition 2.1 mean 
that $\prec$ is {\it two-sided compatible with the multiplication 
operation of the algebra $A$}. Originally in [K-RW], the use of a 
(two-sided) monomial ordering $\prec$  on a solvable polynomial 
algebra $A$ first  guarantees  that {\it $A$ is a domain}, and 
furthermore guarantees an effective (left, right,  two-sided) finite 
Gr\"obner basis theory for $A$ (Theorem 2.3 below). }\v5  

Note that if a $K$-algebra $A=K[a_1,\ldots ,a_n]$ has the PBW
$K$-basis $\B =\{ a^{\alpha}=a_1^{\alpha_1}\cdots
a_n^{\alpha_n}~|~\alpha =(\alpha_1,\ldots ,\alpha_n)\in\NZ^n\}$,
then for any given $n$-tuple $(m_1,\ldots ,m_n)\in\mathbb{N}^n$, a
{\it weighted degree function} $d(~)$ is well defined on nonzero
elements of $A$, namely, for each $a^{\alpha}=a_1^{\alpha_1}\cdots
a_n^{\alpha_n}\in\B$, $d(a^{\alpha})=m_1\alpha_1+\cdots
+m_n\alpha_n$, and for each nonzero
$f=\sum_{i=1}^s\lambda_ia^{\alpha (i)}\in A$ with $\lambda_i\in K^*$
and $a^{\alpha (i)}\in\B$, $d(f)=\max\{ d(a^{\alpha (i)})~|~1\le
i\le s\}.$ If $d(a_i)=m_i>0$ for $1\le i\le n$, then $d(~)$ is
referred to as a {\it positive-degree function} on $A$. \par

Let  $d(~)$ be a  positive-degree function on $A$. If $\prec$ is a
monomial ordering on $\B$ such that for all
$a^{\alpha},a^{\beta}\in\B$,
$$a^{\alpha}\prec a^{\beta}~\hbox{implies}~d(a^{\alpha})\le d(a^{\beta}),\leqno{(*)}$$
then we call $\prec$ a {\it graded monomial ordering} with respect
to $d(~)$, and from now on, unless otherwise stated we always use
$\prec_{gr}$ to denote a  graded monomial ordering.\v5

As one may see from  the literature that in both the commutative and
noncommutative computational algebra, the most popularly used graded
monomial orderings on an algebra $A$ with the PBW $K$-basis $\B$ are
those graded (reverse) lexicographic orderings with respect to the
degree function $d(~)$ such that $d(a_i)=1$, $1\le i\le
n$.\v5

Originally, a noncommutative solvable polynomial algebra (or an 
algebra of solvable type) $R'$ was defined in [K-RW] by first fixing 
a monomial ordering $\prec$ on the standard $K$-basis $\mathscr{B} 
=\{ X_1^{\alpha_1}\cdots X_n^{\alpha_n}~|~\alpha_i\in\NZ\}$ of the 
commutative polynomial algebra  $R=K[X_1,\ldots ,X_n]$ in $n$ 
variables $X_1,\ldots ,X_n$ over a field $K$, and then introducing a 
new multiplication $*$ on $R$, such that certain axioms ([K-RW], 
AXIOMS 1.2) are satisfied. In [LW] the definition of a solvable 
polynomial algebra was modified, in the formal language of 
associative $K$-algebras, as follows. {\parindent=0pt\v5

{\bf 2.2. Definition} Suppose that the $K$-algebra $A=K[a_1,\ldots
,a_n]$ has an admissible system $(\B ,\prec )$. If for all
$a^{\alpha}=a_1^{\alpha_1}\cdots a_n^{\alpha_n}$,
$a^{\beta}=a_1^{\beta_1}\cdots a^{\beta_n}_n\in\B$, the following
condition is satisfied:
$$\begin{array}{rcl} a^{\alpha}a^{\beta}&=&\lambda_{\alpha ,\beta}a^{\alpha +\beta}+f_{\alpha ,\beta},\\
&{~}&\hbox{where}~\lambda_{\alpha ,\beta}\in K^*,~a^{\alpha
+\beta}=a_1^{\alpha_1+\beta_1}\cdots a_n^{\alpha_n+\beta_n},~\hbox{and}\\
&{~}&f_{\alpha ,\beta}\in K\hbox{-span}\B~\hbox{with}~\LM (f_{\alpha
,\beta})\prec a^{\alpha +\beta}~\hbox{whenever}~f_{\alpha ,\beta}\ne
0,\end{array}$$ then $A$ is said to be a {\it solvable polynomial
algebra}. \v5

{\bf Remark} Let $A=K[a_1,\ldots ,a_n]$ be a finitely generated
$K$-algebra and $K\langle X\rangle =K\langle X_1,\ldots ,X_n\rangle$
the free $K$-algebra on $\{X_1,\ldots ,X_n\}$. Then it follows from
[Li4] that $A$ is a solvable polynomial algebra if and only if}\par

(1) $A\cong K\langle X\rangle /\langle G\rangle$ with a finite set
of defining relations $G=\{ g_1,\ldots ,g_m\}$ such that with
respect to some monomial ordering $\prec_{_X}$ on $\KX$, $G$ is a
Gr\"obner basis  of the ideal $\langle G\rangle$ and the set of
normal monomials (mod $G$) gives rise to a PBW $K$-basis $\B$ for
$A$, and\par

(2) there is a monomial ordering $\prec$ on $\B$ such that the
condition on monomials given in  Definition 2.2 is
satisfied.{\parindent=0pt\par

Thus, solvable polynomial algebras are completely determinable and
constructible in a computational way.}\v5

By Definition 2.2 it is straightforward that if $A$ is a solvable
polynomial algebra and $f,g\in A$ with $\LM (f)=a^{\alpha}$, $\LM
(g)=a^{\beta}$, then
$$\LM (fg)=\LM (\LM (f)\LM (g))=\LM (a^{\alpha}a^{\beta})=a^{\alpha +\beta}.\leqno{(\mathbb{P}1)}$$
We shall freely use this property in the rest of this paper without
additional indication.\par

The results mentioned in the theorem below are summarized from
([K-RW], Sections 2 -- 5).{\parindent=0pt\v5

{\bf 2.3. Theorem}  Let $A=K[a_1,\ldots ,a_n]$ be a solvable
polynomial algebra with  admissible system $(\B ,\prec )$. The
following statements hold.\par

(i) $A$ is a (left and right) Noetherian domain.\par

(ii) With respect to the given $\prec$ on $\B$, every nonzero left
ideal $I$ of $A$ has a finite  left Gr\"obner basis $\G=\{
g_1,\ldots ,g_t\}\subset I$ in the sense
that}{\parindent=.5truecm\par

\item{$\bullet$} if $f\in I$ and $f\ne 0$, then there is a $g_i\in\G$ such that $\LM (g_i)|\LM (f)$,
i.e., there is some $a^{\gamma}\in\B$ such that $\LM (f)=\LM
(a^{\gamma}\LM (g_i))$, or equivalently, with $\gamma
(i_j)=(\gamma_{i_{1j}},\gamma_{i_{2j}},\ldots
,\gamma_{i_{nj}})\in\NZ^n$, $f$ has a left Gr\"obner representation:
$$\begin{array}{rcl} f&=&\sum_{i,j}\lambda_{ij}a^{\gamma (i_j)}g_j,~\hbox{where}~\lambda_{ij}\in K^*,
~a^{\gamma (i_j)}\in\B ,~
g_j\in \G,\\
&{~}&\hbox{satisfying}~\LM (a^{\gamma (i_j)}g_j)\preceq\LM
(f)~\hbox{for all}~(i,j).\end{array}$$}{\parindent=0pt\par

(iii) The Buchberger algorithm, that computes a finite Gr\"obner
basis for a finitely generated commutative polynomial ideal, has a
complete noncommutative version that computes a finite left
Gr\"obner basis for a finitely generated left ideal
$I=\sum_{i=1}^mAf_i$ of $A$ (see {\bf Algorithm 1} given in the end
of this section).\par

(iv) Similar results of (ii) and (iii) hold for right ideals and
two-sided ideals of $A$.\par\QED}\v5

Let $A=K[a_1,\ldots ,a_n]$ be a solvable polynomial algebra with
admissible system $(\B ,\prec )$, and let $L=\oplus_{i=1}^sAe_i$ be
a free left $A$-module with the $A$-basis $\{ e_1,\ldots ,e_s\}$.
Then $L$ is a Noetherian module with the $K$-basis $$\BE =\{
a^{\alpha}e_i~|~a^{\alpha}\in\B ,~1\le i\le s\} .$$ We also call
elements of $\BE$ {\it monomials} in $L$. If $\prec_{e}$ is a total 
ordering on $\BE$, and if $\xi =\sum_{j=1}^m\lambda_ja^{\alpha 
(j)}e_{i_j}\in L$, where $\lambda_j\in K^*$ and $\alpha 
(j)=(\alpha_{j_1},\ldots ,\alpha_{j_n})\in\NZ^n$, such that 
$a^{\alpha (1)}e_{i_1}\prec_{e} a^{\alpha 
(2)}e_{i_2}\prec_{e}\cdots\prec_{e} a^{\alpha (m)}e_{i_m}$, then by 
$\LM (\xi )$ we denote the {\it leading monomial} $a^{\alpha 
(m)}e_{i_m}$ of $\xi $, by $\LC (\xi )$ we denote the {\it leading 
coefficient} $\lambda_m$ of $\xi $,  and by $\LT (\xi )$ we denote 
the {\it leading term} $\lambda_ma^{\alpha (m)}e_{i_m}$ of $f$.
\par

With respect to the given monomial ordering $\prec$ on $\B$, a total
ordering $\prec_{e}$ on $\BE$ is called a {\it left monomial 
ordering} if the following two conditions are satisfied:
\par

(1) $a^{\alpha}e_i\prec_{e} a^{\beta}e_j$ implies  $\LM
(a^{\gamma}a^{\alpha}e_i)\prec_{e} \LM (a^{\gamma}a^{\beta}e_j)$ for 
all $a^{\alpha}e_i$, $a^{\beta}e_j\in\BE$, $a^{\gamma}\in\B$;\par

(2) $a^{\beta}\prec a^{\beta}$ implies $a^{\alpha}e_i\prec_{e}  
a^{\beta}e_i$ for all $a^{\alpha},a^{\beta}\in\B$ and $1\le i\le s$. 
{\parindent=0pt\par

From the definition it is straightforward to check that every left
monomial ordering  $\prec_{e}$ on $\BE$ is a well-ordering. 
Moreover, if  $f\in A$ with $\LM (f)=a^{\gamma}$ and $\xi\in L$ with 
$\LM (\xi )=a^{\alpha}e_i$, then by referring to the foregoing 
$(\mathbb{P}1)$ we have
$$\LM (f\xi )=\LM (\LM (f)\LM (\xi ))=\LM (a^{\gamma}a^{\alpha}e_i)=a^{\gamma +\alpha}e_i.\leqno{(\mathbb{P}2)}$$
We shall also freely use this property in the rest of this paper
without additional indication.}\par

Actually as in the commutative case ([AL2], [KR]), any left monomial
ordering $\prec$ on $\B$ may induce  two  left monomial orderings on
$\BE$:
$$\begin{array}{l} (\hbox{{\bf TOP} ordering})\quad a^{\alpha}e_i\prec_{e} a^{\beta}e_j\Leftrightarrow
a^{\alpha}\prec a^{\beta},~\hbox{or}~ a^{\alpha}=a^{\beta}~
\hbox{and}~i<j;\\
(\hbox{{\bf POT} ordering})\quad a^{\alpha}e_i\prec_{e}
a^{\beta}e_j\Leftrightarrow i<j,~
\hbox{or}~i=j~\hbox{and}~a^{\alpha}\prec a^{\beta}.\end{array}$$ \v5

Let $\prec_{e}$ be a left monomial ordering on the $K$-basis $\BE$ 
of $L$, and $a^{\alpha}e_i$, $a^{\beta}e_j\in\BE$, where $\alpha 
=(\alpha_1,\ldots ,\alpha_n)$, $\beta =(\beta_1,\ldots 
,\beta_n)\in\NZ^n$. We say that {\it $a^{\alpha}e_i$ divides 
$a^{\beta}e_j$}, denoted $a^{\alpha}e_i|a^{\beta}e_j$, if $i=j$ and 
$a^{\beta}e_i=\LM (a^{\gamma}a^{\alpha}e_i)$ for some 
$a^{\gamma}\in\B$. It follows from the foregoing property 
$(\mathbb{P}2)$ that
$$a^{\alpha}e_i|a^{\beta}e_j~\hbox{if and only if}~i=j~\hbox{and}~\beta_i\ge
\alpha_i,~1\le i\le n.$$ This division of monomials can be extended
to a division algorithm of dividing an element $\xi$ by a finite
subset of nonzero elements  $U =\{\xi_1,\ldots ,\xi_m\}$ in $L$. 
That is, if there is some $\xi_{i_1}\in U$ such that $\LM 
(\xi_{i_1})|\LM (\xi )$, i.e., there is a monomial $a^{\alpha 
(i_1)}\in\B$ such that $\LM (\xi )=\LM (a^{\alpha (i_1)}\xi_{i_1})$, 
then $\xi ':=\xi -\frac{\LC (\xi )}{\LC (a^{\alpha 
(i_1)}\xi_{i_1})}a^{\alpha (i_1)}\xi_{i_1}$; otherwise, $\xi ':=\xi 
-\LT (\xi )$. Executing this procedure for $\xi '$ and so on, it 
follows from the well-ordering property of $\prec_{e}$ that after 
finitely many repetitions  $\xi$ has an expression
$$\begin{array}{rcl} \xi&=&\sum_{i,j}\lambda_{ij}a^{\alpha (i_j)}\xi_j+\eta ,~\hbox{where}~\lambda_{ij}\in K,
~a^{\alpha (i_j)}\in\B ,~\xi_j\in U ,\\
&{~}&\eta =0~\hbox{or}~\eta =\sum_k\lambda_ka^{\gamma
(k)}e_k~\hbox{with}~\lambda_k\in K,~a^{\gamma (k)}e_k\in\B
(\mathbbm{e})
,\\
&{~}&\hbox{satisfying}\\
&{~}&\LM (a^{\alpha (i_j)}\xi_j)\preceq_{\mathbbm{e}}\LM (\xi
)~\hbox{for
all}~\lambda_{ij}\ne 0,~\hbox{and if}~\eta\ne 0,~\hbox{then}\\
&{~}&a^{\gamma (k)}e_k\preceq_{\mathbbm{e}}\LM (\xi),~\LM
(\xi_i){\not |}~a^{\gamma (k)}e_k~\hbox{for all}~\xi_i\in 
U~\hbox{and all}~\lambda_k\ne 0.\end{array}$$ The element $\eta$ 
appeared in the above expression is called a {\it remainder} of 
$\xi$ on division by $ U$, and is usually denoted by $\OV{\xi}^{ 
U}$, i.e., $\OV{\xi}^{ U}=\eta$. If $\OV{\xi}^{ U}=0$, then we say 
that $\xi$ is {\it reduced to zero}  on division by $ U$. A nonzero 
element $\xi\in L$ is said to be {\it normal} (mod $ U$) if $\xi 
=\OV{\xi}^{ U}$.\par

Based on the division algorithm, the notion of a {\it left Gr\"obner
basis} for a submodule $N$ of the free  module
$L=\oplus_{i=1}^sAe_i$ comes into play. Since $A$ is a Noetherian
domain, it follows that  $L$ is a Noetherian $A$-module and the
following proposition holds. {\parindent=0pt\v5

{\bf 2.4. Theorem} With respect to the given $\prec_{e}$ on $\BE$, 
every nonzero submodule $N$ of $L$ has a finite left Gr\"obner basis 
$\G =\{ g_1,\ldots ,g_m\}\subset N$ in the sense that 
{\parindent=.5truecm\par

\item{$\bullet$} if $\xi\in N$ and $\xi\ne 0$,
then $\LM (g_i)|\LM (\xi )$ for some $g_i\in\G$, i.e., there is a
monomial $a^{\gamma}\in\B$ such that $\LM (\xi )=\LM (a^{\gamma}\LM
(g_i))$, or equivalently, $\xi$ has a {\it left Gr\"obner
representation} $\xi =\sum_{i,j}\lambda_{ij}a^{\alpha (i_j)}g_j$,
where $\lambda_{ij}\in K^*$, $a^{\alpha (i_j)}\in\B$ with $\alpha
(i_j)=(\alpha_{i_{j1}},\ldots ,\alpha_{i_{jn}})\in\NZ^n$,
$g_j\in\G$,  satisfying $\LM (a^{\alpha
(i_j)}g_j)\preceq_{\mathbbm{e}}\LM (\xi)$;\par}

moreover, starting with any finite generating set of $N$, such a
left Gr\"obner basis $\G$ can be computed by running a
noncommutative version of the Buchberger algorithm for modules over
solvable polynomial algebras (see {\bf Algorithm 1} presented below). \par\QED}\v5  

Since the noncommutative version of Buchberger algorithm is based on
the noncommutative version of Buchberger's  criterion that makes the
strategy for computing left Gr\"obner  bases of modules over
solvable polynomial algebras, for the reader's convenience and the
use of the next section, we recall both of them  as follows.\par

Let $N=\sum_{i=1}^mA\xi_i$ with $ U =\{\xi_1,\ldots ,\xi_m\}\subset 
L$.  For  $\xi_i,\xi_j\in U$ with $1\le i<j\le m$, $\LM 
(\xi_i)=a^{\alpha (i)}e_{p}$, $\LM (\xi_j)=a^{\alpha (j)}e_{q}$, 
where $\alpha (i)=(\alpha_{i_1},\ldots,\alpha_{i_n})$, $\alpha 
(j)=(\alpha_{j_1},\ldots ,\alpha_{j_n})$,  put $\gamma 
=(\gamma_1,\ldots ,\gamma_n)$ with $\gamma_k =\max\{ 
\alpha_{i_k},\alpha_{j_k})$. The {\it left S-polynomial} of $\xi_i$ 
and $\xi_j$ is defined as
$$S_{\ell}(\xi_i,\xi_j)=\left\{\begin{array}{ll}
\displaystyle{\frac{1}{\LC (a^{\gamma -\alpha (i)}\xi_i)}}a^{\gamma
-\alpha (i)}\xi_i-\displaystyle{\frac{1}{\LC (a^{\gamma -\alpha
(j)}\xi_j)}}a^{\gamma
-\alpha (j)}\xi_j,&\hbox{if}~ p=q\\
0,&\hbox{if}~p\ne q.\end{array}\right.$$ {\parindent=0pt\vskip
6pt\def\S{{\cal S}}  

{\bf 2.5. Theorem} (Noncommutative version of Buchberger's
criterion) With notation as above, $U$ is a left Gr\"obner basis of
the submodule $N$ if and only if every $S_{\ell}(\xi_i,\xi_j)$ is
reduced to 0 on division by $U$, i.e.,
$\OV{S_{\ell}(\xi_i,\xi_j)}^U=0$. \vskip 6pt

\underline{{\bf Algorithm 1} (Noncommutative version of Buchberger 
algorithm)~~~~~~~~~~~~~~~~~~~~~~~~~~~~~~~~~~~~~~~~~~~~}\par

$\begin{array}{l} \textsc{INPUT}: ~ U =
\{ \xi_1,...,\xi_m\}\\
\textsc{OUTPUT}:~ \G =\{ g_1,...,g_t\}~\hbox{a left Gr\"obner basis
of}~N=\sum_{i=1}^mA\xi_i~\\
\textsc{INITIALIZATION}:~ m':=m,~\G :=\{g_1=\xi_1,\ldots ,g_{m'}=\xi_m\}  ,\\
\hskip 4.22truecm \S :=\left\{S_{\ell}(g_i,g_j)~\left
|~\begin{array}{l} g_i,g_j\in\G ,~i<j,~\hbox{and for some}~e_t,\\
\LM (g_i)=a^{\alpha}e_t,~\LM
(g_j)=a^{\beta}e_t\end{array}\right.\right\}\end{array}$

$\begin{array}{l}
\textsc{BEGIN}\\
~~~~\textsc{WHILE}~\S\ne \emptyset\\
~~~~~~~~~~\hbox{Choose any}~S_{\ell}(g_i,g_j)\in\S\\
~~~~~~~~~~\S :=\S -\{ S_{\ell}(g_i,g_j)\}\\
~~~~~~~~~~\OV{S_{\ell}(g_i,g_j)}^{\G}=\eta\\
~~~~~~~~~~~~~\textsc{IF}~\eta\ne 0~\hbox{with}~\LM (\eta )=
a^{\rho}e_k~\textsc{THEN}\\
~~~~~~~~~~~~~~~~~m':=m'+1,~g_{m'}:=\eta \\
~~~~~~~~~~~~~~~~~\S :=\S\cup\{ S_{\ell}(g_j,g_{m'})~|~g_j\in\G,~\LM (g_j )=a^{\nu}e_k\}\\
~~~~~~~~~~~~~~~~~\G :=\G\cup\{g_{m'}\},\\
~~~~~~~~~~~~~\textsc{END}\\
~~~~\textsc{END}\\
\textsc{END}\end{array}$\par
\underline{~~~~~~~~~~~~~~~~~~~~~~~~~~~~~~~~~~~~~~~~~~~~~~~~~~~~~~~~~~~~~~~~~~~~~~~~~~~~~~~
~~~~~~~~~~~~~~~~~~~~~~~~~~~~~~~~~~~~~~~~~~~~~~~~~~~} \vskip 6pt

One is referred to the up-to-date computer algebra systems
\textsc{Singular} [DGPS] for the implementation of {\bf Algorithm
1}. Also, nowadays there have been optimized algorithms, such as the
signature-based algorithm for computing Gr\"obner bases in solvable
polynomial algebras [SWMZ], which is based on the celebrated F5
algorithm [Fau],  may be used to speed-up the computation of left 
Gr\"obner bases for modules.}  \v5

\section*{3. Computation of Minimal Homogeneous Generating Sets of Graded Submodules $N\subset L$}
Based on {\bf Algorithm  1} presented in the last section, in this
section we  show that the methods and algorithms, developed in
([CDNR], [KR]) for computing minimal homogeneous generating sets of
graded submodules in free modules over commutative polynomial
algebras, can be adapted for  computing minimal homogeneous
generating sets of graded submodules in free modules over weighted
$\NZ$-graded solvable polynomial algebras. All notions, notations
and conventions introduced before are maintained. \v5

Let $A=K[a_1,\ldots ,a_n]$ be a solvable polynomial $K$-algebra with
admissible system $(\B ,\prec )$, where $\B =\{
a^{\alpha}=a_1^{\alpha_1}\cdots a_n^{\alpha_n}~|~\alpha
=(\alpha_1,\ldots ,\alpha_n)\in\NZ^n\}$ is the PBW $K$-basis of $A$
and $\prec$ is a monomial ordering on $\B$. Given a positive-degree
function $d(~)$ on $A$ (as defined in  Section 2) such that
$d(a_i)=m_i>0$, $1\le i\le n$, we know that $A$ has  an $\NZ$-graded
$K$-module structure, i.e., $A=\oplus_{p\in\NZ}A_p$ with
$A_p=K\hbox{-span}\{ a^{\alpha}\in\B~|~d(a^{\alpha})=p\}$, in
particular, $A_0=K$; if furthermore $A_{p_1}A_{p_2}\subseteq
A_{p_1+p_2}~\hbox{for all}~p_1,p_2\in\NZ ,$ then $A$ is turned into
a weighted {\it $\NZ$-graded solvable polynomial algebra} with
respect to $d(~)$. In this case, elements in $A_p$ are called {\it
homogeneous elements of degree $p$}, and accordingly $A_p$ is called
the {\it degree-$p$ homogeneous part} of $A$. For every nonzero
$h\in A_p$, we  write $d_{\rm gr}(h)$ for the degree of $h$ as a
homogeneous element of $A$, i.e., $d_{\rm gr}(h)=p$. Since  $A$ is a
domain, the degree function $d_{\rm gr}(~)$ defined on nonzero
homogeneous elements of $A$  has the property:
$$d_{\rm gr}(h_1h_2)=d_{\rm gr}(h_1)+d_{\rm gr}(h_2)~\hbox{for all nonzero homogeneous elements}
~h_1,h_2\in A.\leqno{(\mathbb{P}3)}$$ From now on we shall freely use this
property without additional indication.{\parindent=0pt \v5

{\bf Remark} (i) We emphasize that if $A=K[a_1,\ldots ,a_n]$ is a
weighted $\NZ$-graded solvable polynomial algebra with respect to a
positive degree function $d(~)$ on $A$, then {\it every
$a^{\alpha}\in\B$ is a homogeneous elements of $A$ and
$d(a^{\alpha})=d_{\rm gr}(a^{\alpha})$}.\par

(ii) Given a solvable polynomial algebra $A=K[a_1,\ldots ,a_n]$ with
admissible system $(\B ,\prec )$,  it follows from Definition 2.2
that{\parindent=1truecm \vskip6pt

\item{$\bullet$} $A$ is a weighted $\NZ$-graded solvable polynomial algebra with respect to a given
positive-degree function $d(~)$ if and only if for   $1\le i<j\le
n$, in the relation   $a_ja_i=\lambda_{ji}a_ia_j+f_{ji}$ with
$f_{ji}=\sum\mu_ka^{\alpha (k)}$, $d(a^{\alpha (k)})=d(a_ia_j)$
whenever $\mu_k\ne 0$. \par}{\parindent=0pt\vskip 6pt

Consequently, if $\prec_{gr}$ is a graded monomial ordering on $\B$
with respect to some given positive-degree function $d(~)$ on $A$,
then it is easy to know whether $A$ is a weighted $\NZ$-graded
algebra with respect to $d(~)$ or not.}}\v5

In view of the above remark, typical $\NZ$-graded solvable
polynomial algebras  are those iterated skew polynomial $K$-algebras
$A=K[a_1,\ldots ,a_n]$ subject to the relations
$$a_ja_i=\lambda_{ji}a_ia_j,~\quad\lambda_{ji}\in K^*,~1\le i<j\le n,$$
where the positive-degree function on $A$ can be defined by
$d(a_i)=m_i$ for any fixed tuple $(m_1,\ldots ,m_n)$ of positive
integers.  Such algebras include the well-known coordinate rings of
quantum affine $n$-spaces. Another well-known $\NZ$-graded solvable
polynomial algebra is the coordinate ring $M_q(2)=K[a,b,c,d]$ of the
manifold
of quantum $2\times 2$ matrices $\left (\begin{array}{cc} a&b\\
c&d\end{array}\right )$ introduced in [Man], which has the defining
relations
$$\begin{array}{ll} ab=q^{-1}ba,&db=bd-(q-q^{-1})ac,\\
cb=qbc,&da=qad,\\
ca=ac,&dc=qcd,\end{array}$$ where each generator is assigned the
degree 1. More generally, let $\Lambda =(\lambda_{ij})$ be a
multiplicatively antisymmetric $n\times n$ matrix over $K$, and let
$\lambda\in K^*$ with $\lambda\ne -1$. Considering the
multiparameter coordinate ring of quantum $n\times n$ matrices over
$K$ (see [Good]), namely the $K$-algebra ${\cal O}_{\lambda
,\Lambda}(M_n(K))$ generated by $n^2$ elements $a_{ij}$ ($1\le
i,j\le n$) subject to the relations
$$a_{\ell m}a_{ij}=\left\{\begin{array}{ll} \lambda_{\ell i}\lambda_{jm}a_{ij}a_{\ell m}+
(\lambda -1)\lambda_{\ell i}a_{im}a_{\ell j}&(\ell >i,~m>j)\\
\lambda\lambda_{\ell i}\lambda_{jm}a_{ij}a_{\ell m}&(\ell >i,~m\le
j)\\ \lambda_{jm}a_{ij}a_{\ell m}&(\ell =i,~m>j)\end{array}\right.$$
Then ${\cal O}_{\lambda ,\Lambda}(M_n(K))$ is an $\NZ$-graded
solvable polynomial algebra, where each generator has degree 1. The
two examples given below provide weighted $\NZ$-graded solvable
polynomial algebras in which some generators may have degree $\ge
2$. Moreover, it is also known that the associated graded algebra
and the Rees algebra of every $\NZ$-filtered solvable polynomial
algebra with a graded monomial ordering are $\NZ$-graded solvable
polynomial algebras (see [LW], [Li1]).\v5

Let $A=K[a_1,\ldots ,a_n]$ be a weighted $\NZ$-graded solvable
polynomial algebra with respect to a given positive-degree function
$d(~)$ on $A$, and let $(\B ,\prec )$ be an admissible system of
$A$. If $L=\oplus_{i=1}^sAe_i$ is a free left $A$-module with the
$A$-basis $\{ e_1,\ldots ,e_s\}$, then $L$ has the $K$-basis $\B
(\mathbbm{e})=\{ a^{\alpha}e_i~|~a^{\alpha}\in\B ,~1\le i\le n\}$,
and for an {\it arbitrarily} fixed $\{b_1,\ldots ,b_s\}\subset\NZ$,
$L$ can be turned into an {\it $\NZ$-graded free $A$-module} in the
sense that $L=\oplus_{q\in\NZ}L_q$ in which
$$L_q=\{ 0\}~\hbox{if}~q<\min\{ b_1,\ldots ,b_s\};~\hbox{otherwise}~L_q=
\sum_{p_i+b_i=q}A_{p_i}e_i,\quad q\in\NZ ,$$
or alternatively, for $q\ge \min\{ b_1,\ldots ,b_s\}$,
$$L_q=K\hbox{-span}\{ a^{\alpha}e_i\in\BE~|~d(a^{\alpha})+b_i=q\},
~q\in \NZ ,$$
such that $A_pL_q\subseteq A_{p+q}$ for all $p,q\in\NZ$. For each
$q\in\NZ$, elements in $L_q$ are called {\it homogeneous elements of
degree $q$}, and accordingly $L_q$ is called the {\it degree-$q$
homogeneous part} of $L$.  If $\xi\in L_q$, then we write $d_{\rm
gr}(\xi )$ for the degree of $\eta$ as a homogeneous element of $L$,
i.e., $d_{\rm gr}(\xi )=q$. In particular, $d_{\rm gr}(e_i)=b_i$,
$1\le i\le s$. As with the degree of homogeneous elements in $A$,
noticing that $d_{\rm gr}(a^{\alpha}e_i)=d(a^{\alpha})+b_i$ for all
$a^{\alpha}e_i\in\BE$ and that $A$ is a domain, from now on we shall
freely use the following property without additional indication: for
all nonzero homogeneous elements $h\in A$ and all nonzero
homogeneous elements $\xi\in L$,
$$d_{\rm gr}(h\xi )=d_{\rm gr}(h)+d_{\rm gr}(\xi ).\leqno{(\mathbb{P}4)}$${\parindent=0pt\par

{\bf Convention} Unless otherwise stated, from now on throughout the
subsequent sections if we say that $L$ is a graded free module over
a weighted $\NZ$-graded solvable polynomial algebra $A$, then it
always means that  $L$ has the $\NZ$-gradation as constructed
above.}\v5

Let $L=\oplus_{q\in\NZ}L_q$ be a graded free $A$-module. If a
submodule $N$ is generated by homogeneous elements, then $N$ is
called a {\it graded submodule} of $L$. A graded submodule $N$ has
the $\NZ$-graded structure $N=\oplus_{q\in\NZ}N_q$ with $N_q=N\cap
L_q$, such that $A_pN_q\subseteq N_{p+q}$ for all $p,q\in\NZ$. Note
that monomials in $\B$ are homogeneous elements of $A$, thereby left
$S$-polynomials of homogeneous elements are homogeneous elements,
and remainders of homogeneous elements on division by homogeneous
remain homogeneous elements. It follows that if a graded submodule
$N=\sum_{i=1}^mA\xi_i$ of $L$ is generated by the set of homogeneous
elements $\{ \xi_1,\ldots ,\xi_m\}$, then, running  the
noncommutative version of Buchberger's algorithm for modules over
solvable polynomial algebras ({\bf Algorithm 1} in Section 2) with
respect to a fixed monomial ordering $\prec_{e}$ on $\B
(\mathbbm{e})$, it produces a finite {\it homogeneous left Gr\"obner
basis} $\G$ for $N$, that is, $\G$ consists of homogeneous
elements.\v5

In what follows, $A=K[a_1,\ldots ,a_n]$ denotes a weighted
$\NZ$-graded solvable polynomial algebra with respect to a given
positive-degree function $d(~)$ on $A$, $(\B,\prec )$ denotes a
fixed admissible system of $A$, $L=\oplus_{i=1}^sAe_i$ denotes an
$\NZ$-graded free $A$-module such that $d_{\rm gr}(e_i)=b_i$, $1\le
i\le s$, and  $\prec_{e}$ denotes a fixed left monomial ordering on 
the $K$-basis $\BE$ of $L$. Moreover,  as in Section 2 we write 
$S_{\ell}(\xi_i,\xi_j)$ for the left S-polynomial of two elements 
$\xi_i$, $\xi_j\in L$.\par

To reach our goal of this section on the basis of {\bf Algorithm 1},
let us point out that although  monomials from the PBW $K$-basis
$\B$ of a weighted $\NZ$-graded solvable polynomial algebra $A$ {\it
can no longer behave} as well as monomials in a commutative
polynomial algebra (namely the product of two monomials is not
necessarily a monomial), every monomial from $\B$ is a homogeneous
element in the weighted $\NZ$-graded structure of $A$ (as we
remarked before), thereby the product of two monomials is a
homogeneous element. Bearing in mind this fact, the rule of
division, and the properties $(\mathbb{P}1) - (\mathbb{P}3)$
mentioned in Section 2 and the foregoing  $(\mathbb{P}4)$, the
argument below will go through without trouble.\par

We first show that the algorithm given in ([KR], Proposition 4.5.10)
can be modified to compute $n$-truncated left Gr\"obner bases for
graded submodules of free modules over a weighted $\NZ$-graded
solvable polynomial algebra. The discussion on $n$-truncated left
Gr\"obner bases presented below is similar to the commutative case
as in [KR]. \def\S{{\cal S}}{\parindent=0pt\v5

{\bf 3.1. Definition} Let $G=\{ g_1,\ldots ,g_t\}$ be a subset of
homogeneous elements of $L$, $N=\sum_{i=1}^tAg_i$ the graded
submodule generated by $G$, and let $n\in\NZ$, $G_{\le n}=\{ g_j\in
G~|~d_{\rm gr}(g_j)\le n\}$. If, for each nonzero homogeneous
element $\xi\in N$ with $d_{\rm gr}(\xi )\le n$, there is some
$g_i\in G_{\le n}$ such that $\LM (g_i)|\LM (\xi )$ with respect to
$\prec_{e}$, then we call $G_{\le n}$ an $n$-{\it truncated left 
Gr\"obner basis} of $N$ with respect to $(\BE, \prec_{e})$.}\v5

Verification of the lemma below is straightforward.
{\parindent=0pt\v5

{\bf 3.2. Lemma} Let $\G =\{ g_1,\ldots ,g_t\}$ be a homogeneous
left Gr\"obner basis for the graded submodule $N=\sum_{i=1}^tAg_i$
of $L$ with respect to $(\BE ,\prec_{e})$. For each $n\in\NZ$, put 
$\G_{\le n}=\{ g_j\in\G~|~d_{\rm gr}(g_j)\le n\}$, $N_{\le 
n}=\cup_{q=0}^nN_q$ where each $N_q$ is the degree-$q$ homogeneous 
part of $N$, and let $N(n)=\sum_{\xi\in N_{\le n}}A\xi$ be  the 
graded submodule generated by $N_{\le n}$.  The following statements 
hold.\par

(i) $\G_{\le n}$ is an $n$-truncated left Gr\"obner basis of $N$.
Thus, if $\xi\in L$ is a homogeneous element with $d_{\rm gr}(\xi
)\le n$,  then $\xi\in N$ if and only if $\OV{\xi}^{\G_{\le n}}=0$,
i.e., $\xi$ is reduced to zero on division by $\G_{\le n}$,.\par

(ii) $N(n)=\sum_{g_j\in\G_{\le n}}Ag_j$, and $\G_{\le n}$ is an
$n$-truncated left Gr\"obner basis of $N(n)$.\par\QED}\v5

In light of {\bf Algorithm 1} (presented in Section 2), an
$n$-truncated left Gr\"obner basis is characterized as
follows.{\parindent=0pt\v5

{\bf 3.3. Proposition} Let $N=\sum_{i=0}^sAg_i$ be the graded
submodule of $L$ generated by a set of homogeneous elements $G=\{
g_1,\ldots ,g_m\}$. For each $n\in\NZ$, put  $G_{\le n}=\{ g_j\in
G~|~d_{\rm gr}(g_j)\le n\}$. The following statements are equivalent
with respect to the given $(\BE, \prec_{e})$.\par

(i) $G_{\le n}$ is an $n$-truncated left Gr\"obner basis of $N$.\par

(ii) Every nonzero left S-polynomial $S_{\ell}(g_i,g_j)$ of $d_{\rm
gr}(S_{\ell}(g_i,g_j))\le n$ is reduced to zero on division by
$G_{\le n}$, i.e, $\OV{S_{\ell}(g_i,g_j)}^{G_{\le n}}=0$.\vskip 6pt

{\bf Proof} Recall that if $g_i,g_j\in G$,  $\LT
(g_i)=\lambda_ia^{\alpha}e_t$ with $\alpha =(\alpha_1,\ldots
,\alpha_n)$,  $\LT (g_j)=\lambda_ja^{\beta}e_t$ with $\beta
=(\beta_1,\ldots ,\beta_n)$, and $\gamma =(\gamma _1,\ldots
\gamma_n)$ with $\gamma _i=\max\{\alpha_i,\beta_i\}$, $1\le i\le n$,
then
$$S_{\ell}(g_i,g_j)=\frac{1}{\LC(a^{\gamma -\alpha}g_i)}a^{\gamma -\alpha}g_i-
\frac{1}{\LC(a^{\gamma -\beta}g_j)}a^{\gamma -\beta}g_j$$ is a
homogeneous element in $N$ with $d_{\rm
gr}(S_{\ell}(g_i,g_j))=d(a^{\gamma})+b_t$ by the foregoing property
($\mathbb{P}4$). If $d_{\rm gr}(S_{\ell}(g_i,g_j))\le n$, then it
follows from (i) that (ii) holds.}
\par

Conversely, suppose that (ii) holds. To see that $G_{\le n}$ is an
$n$-truncated left Gr\"obner basis of $N$, let us run {\bf Algorithm
1} with the initial input data $G$. Without optimizing {\bf
Algorithm 1} we may certainly assume that $G\subseteq\G$, thereby
$G_{\le n}\subseteq\G_{\le n}$ where $\G$ is the new input set
returned by each pass through the WHILE loop. On the other hand, by
the construction of $S_{\ell}(g_i,g_j)$ and the foregoing property
($\mathbb{P}4$) we know that if $d_{\rm gr}(S_{\ell}(g_i,g_j))\le
n$, then $d_{\rm gr}(g_i)\le n$, $d_{\rm gr}(g_j)\le n$. Hence, the
assumption (ii) implies that {\bf Algorithm 1} does not append any
new element of degree $\le n$ to $\G$. Therefore, $G_{\le n}=\G_{\le
n}$. By Lemma 3.2 we conclude that $G_{\le n}$ is an $n$-truncated
left Gr\"obner basis of $N$.\QED {\parindent=0pt\v5

{\bf 3.4. Corollary} Let  $N=\sum^m_{i=1}Ag_i$ be the graded
submodule of $L$ generated by a set of homogeneous elements $G=\{
g_1,\ldots ,g_m\}$. Suppose that $G_{\le n}=\{ g_j\in G~|~d_{\rm
gr}(g_j)\le n\}$ is an $n$-truncated left Gr\"obner basis of $N$
with respect to $(\BE, \prec_{e})$. \par

(i) If $\xi\in L$ is a nonzero homogeneous element of $d_{\rm
gr}(\xi )=n$ such that $\LM (g_i){\not |}~\LM (\xi )$ for all
$g_i\in G_{\le n}$, then $G'=G_{\le n}\cup\{\xi\}$ is an
$n$-truncated left Gr\"obner basis for both the graded submodules
$N'=N+A\xi$ and $N''=\sum_{g_j\in G_{\le n}}Ag_j+A\xi$ of $L$.\par

(ii) If $n\le n_1$ and $\xi\in L$ is a nonzero homogeneous element
of $d_{\rm gr}(\xi )=n_1$ such that $\LM (g_i){\not |}~\LM (\xi )$
for all $g_i\in G_{\le n}$, then $G'=G_{\le n}\cup\{\xi\}$ is an
$n_1$-truncated left Gr\"obner basis for the graded submodule
$N'=\sum_{g_j\in G_{\le n}}Ag_j+A\xi$ of $L$. \vskip 6pt

{\bf Proof} If $\xi\in L$ is a nonzero homogeneous element of
$d_{\rm gr}(\xi )=n_1\ge n$ and $\LM (\xi_i){\not |}~\LM (\xi)$ for
all $\xi_i\in G_{\le n}$, then noticing the property $(\mathbb{P}2)$
mentioned in Section 2 and the foregoing $(\mathbb{P}4)$, we see
that every nonzero left S-polynomial $S_{\ell}(\xi ,\xi_i)$ with
$\xi_i\in G$ has $d_{\rm gr}(S_{\ell}(\xi ,\xi_i))>n$. Hence both
(i) and (ii) hold by Proposition 3.3.\QED \v5

{\bf 3.5. Proposition} (Compare with ([KR2], Proposition 4.5.10)) 
Given a finite set of nonzero homogeneous elements $ U =\{ 
\xi_1,\ldots ,\xi_m\}\subset L$ with $d_{\rm gr}(\xi_1)\le d_{\rm 
gr}(\xi_2)\le\cdots d_{\rm gr}(\xi_m)$, and a positive integer 
$n_0\ge  d_{\rm gr}(\xi_1)$, the following algorithm computes an 
$n_0$-truncated left Gr\"obner basis $\G =\{ g_1,...,g_t\}$ for the 
graded submodule $N=\sum^m_{i=1}A\xi_i$ such that $d_{\rm 
gr}(g_1)\le d_{\rm gr}(g_2)\le\cdots d_{\rm gr}(g_t)$. 
{\parindent=0pt\vskip 6pt

\underline{\bf Algorithm 2 
~~~~~~~~~~~~~~~~~~~~~~~~~~~~~~~~~~~~~~~~~~~~~~~~~~~~~~~~~~~~~~~~~~~~~~~~~~~~~~~~~~~~~~~~~~~~~~~~~}\par 
$\begin{array}{l} \textsc{INPUT}:~  U = \{ 
\xi_1,...,\xi_m\}~\hbox{with}~d_{\rm gr}(\xi_1)\le d_{\rm 
gr}(\xi_2)\le\cdots d_{\rm
gr}(\xi_m),\\
~~~~~~~~~~~~~~~n_0,~\hbox{where}~n_0\ge d_{\rm gr}(\xi_1)\\
\textsc{OUTPUT}: ~\G =\{ g_1,...,g_t\}~\hbox{an}~n_0\hbox{-truncated 
left Gr\"obner basis
of}~N\\
\textsc{INITIALIZATION}:~ \S_{\le n_0} :=\emptyset ,~W := U,~\G 
:=\emptyset
,~t':=0\end{array}$ \\
$\begin{array}{l}
\textsc{LOOP}\\
n:=\min\{d_{\rm gr}(\xi_i),~d_{\rm gr}(S_{\ell}(g_i,g_j))~|~\xi_i\in 
W
,~S_{\ell}(g_i,g_j)\in\S_{\le n_0}\}\\
\S_n:=\{ S_{\ell}(g_i,g_j)\in\S_{\le n_0}~|~d_{\rm gr}(S_{\ell}(g_i,g_j))=n\},~W_n:=\{\xi_j\in W~|~d_{\rm gr}(\xi_j)=n\}\\
\S_{\le n_0} :=\S_{\le n_0} -\S_n,~W:=W-W_n\\
~~~~~~\textsc{WHILE}~\S_n\ne\emptyset~\textsc{DO}\\
~~~~~~~~~~~~~\hbox{Choose any}~S_{\ell}(g_i,g_j)\in\S_{n}\\
~~~~~~~~~~~~~\S_n :=\S_n -\{ S_{\ell}(g_i,g_j)\}\\
~~~~~~~~~~~~~\OV{S_{\ell}(g_i,g_j)}^{\G}=\eta\\
~~~~~~~~~~~~~\textsc{IF}~\eta\ne 0~\hbox{with}~\LM (\eta )=
a^{\rho}e_k\textsc{THEN}\\
~~~~~~~~~~~~~~~~~~t':=t'+1,~g_{t'}:=\eta\\
~~~~~~~~~~~~~~~~~~\S_{\le n_0} :=\S_{\le n_0}\cup\left\{ 
S_{\ell}(g_i,g_{t'})\left |\begin{array}{l} g_i\in\G ,1\le i< t',\LM
(g_i )=a^{\tau}e_k,\\
d_{\rm gr}(S_{\ell}(g_i,g_{t'}))\le n_0\end{array}\right.\right\}\\
~~~~~~~~~~~~~~~~~~\G :=\G\cup\{ g_{t'}\}\\
~~~~~~~~~~~~~\textsc{END}\\
~~~~~~\textsc{END}\\
\end{array}$

$\begin{array}{l}
~~~~~~\textsc{WHILE}~W_n\ne\emptyset~\textsc{DO}\\
~~~~~~~~~~~~~\hbox{Choose any}~\xi_j\in W_{n}\\
~~~~~~~~~~~~~W_n :=W_n -\{ \xi_j\}\\
~~~~~~~~~~~~~\OV{\xi_j}^{\G}=\eta\\
~~~~~~~~~~~~~\textsc{IF}~\eta\ne 0~\hbox{with}~\LM (\eta )=
a^{\rho}e_k~\textsc{THEN}\\
~~~~~~~~~~~~~~~~~~ t':=t'+1,~g_{t'}:=\eta\\
~~~~~~~~~~~~~~~~~~\S_{\le n_0} :=\S_{\le n_0}\cup\left\{ 
S_{\ell}(g_i,g_{t'})\left |\begin{array}{l} g_i\in\G ,1\le i<t',\LM
(g_i )=a^{\tau}e_k,\\
d_{\rm gr}(S_{\ell}(g_i,g_{t'}))\le n_0\end{array}\right.\right\}\\
~~~~~~~~~~~~~~~~~~\G :=\G\cup\{ g_{t'}\}\\
~~~~~~~~~~~~~\textsc{END}\\
~~~~~~\textsc{END}\\
\textsc{UNTIL}~ \S_{\le n_0} =\emptyset\\
\textsc{END}
\end{array}$
\par
\underline{~~~~~~~~~~~~~~~~~~~~~~~~~~~~~~~~~~~~~~~~~~~~~~~~~~~~~~~~~~~~~~~~~~~~~~~~~~~~~~~~~~~~~~~~~~~~~~~~~~~~~~~~~~~~~~~~~~~~~~~~~~~~~~~~~~~}
} \vskip 6pt

{\bf Proof} First note that both the WHILE loops append new elements
to $\G$ by taking the nonzero normal remainders on division by $\G$.
Thus, with a fixed $n$, by the definition of a left S-polynomial and
the normality of $g_{t'}$ (mod $\G$), it is straightforward to check
that in both the WHILE loops every newly appended
$S_{\ell}(g_i,g_{t'})$ has $d_{\rm gr}(S_{\ell}(g_i,g_{t'}))>n$. To
proceed, let us write $N(n)$ for the submodule generated by $\G$
which is obtained after $W_n$ is exhausted in the second WHILE loop.
If $n_1$ is the first number after $n$ such that
$\S_{n_1}\ne\emptyset$, and for some $S_{\ell}(g_i,g_j)\in
\S_{n_1}$, $\eta =\OV{S_{\ell}(g_i,g_j)}^{\G}\ne 0$ in a certain
pass through the first WHILE loop, then we note that this $\eta$ is
still contained in $N(n)$. Hence, after $\S_{n_1}$ is exhausted in
the first WHILE loop, the obtained $\G$ generates $N(n)$ and $\G$ is
an $n_1$-truncated left Gr\"obner basis of $N(n)$. Noticing that the
algorithm starts with $\S =\emptyset$ and $\G =\emptyset$,
inductively it follows from Proposition 3.3 and Corollary 3.4 that
after $W_{n_1}$ is exhausted in the second WHILE loop, the obtained
$\G$ is an $n_1$-truncated left Gr\"obner basis of $N(n_1)$. Since
$n_0$ is finite and all the generators of $N$ with  $d_{\rm
gr}(\xi_j)\le n_0$ are processed through the second WHILE loop, the
algorithm terminates and the eventually obtained $\G$ is an
$n_0$-truncated left Gr\"obner basis of $N$. Finally, the fact that
the degrees of elements in $\G$ are non-decreasingly ordered follows
from the choice of the next $n$ in the algorithm. \QED }\v5

Let the data $(A,\B,\prec )$ and $(L, \B(\mathbbm{e}),\prec_{e})$ be 
as fixed before. Combining the foregoing results, we now proceed to 
show that the algorithm given in ([KR], Theorem 4.6.3)) can be 
modified to compute  minimal homogeneous generating sets of graded 
submodules in free modules over $A$.\par

Let $N$ be a graded submodule of the $\NZ$-graded free $A$-module
$L$ fixed above. We say that a homogeneous generating set $ U$ of
$N$ is a {\it minimal homogeneous generating set} if any proper
subset of $ U$ cannot be a generating set of $N$. As preparatory
result, we first show that the noncommutative analogue of ([KR],
Proposition 4.6.1, Corollary 4.6.2) holds true for
$N$.{\parindent=0pt\v5

{\bf 3.6. Proposition} Let $N=\sum^m_{i=1}A\xi_i$ be the graded
submodule of $L$ generated by a set of homogeneous elements $ U=\{ 
\xi_1,\ldots ,\xi_m\}$, where $d_{\rm gr}(\xi_1)\le d_{\rm 
gr}(\xi_2)\le \cdots\le d_{\rm gr}(\xi_m)$. Put $N_1=\{ 0\}$, 
$N_i=\sum^{i-1}_{j=1}A\xi_j$, $2\le i\le m$. The following 
statements hold.\par

(i) $ U$ is a minimal homogeneous generating set of $N$ if and only 
if $\xi_i\not\in N_i$, $1\le i\le m$.\par

(ii) The set $\OV U=\{ \xi_k~|~\xi_k\in U ,~\xi_k\not\in N_k\}$ is a 
minimal homogeneous generating set of $N$.\vskip 6pt

{\bf Proof} (i) If $ U$ is a minimal homogeneous generating set of 
$N$, then clearly $\xi_i\not\in N_i$, $1\le i\le m$.}\par

Conversely, suppose $\xi_i\not\in N_i$, $1\le i\le m$. If $ U$ is
not a minimal homogeneous generating set of $N$, then, there is some
$i$ such that $N$ is generated by $\{ \xi_1,\ldots
,\xi_{i-1},\xi_{i+1},\ldots ,\xi_m\}$, thereby $\xi_i=\sum_{j\ne
i}h_j\xi_j$ for some nonzero homogeneous elements $h_j\in A$ such
that $d_{\rm gr}(\xi_i)=d_{\rm gr}(h_jg_j)=d_{\rm gr}(h_j)+d_{\rm
gr}(\xi_j)$, where the second equality follows from the foregoing
property $(\mathbb{P}4)$. Thus $d_{\rm gr}(\xi_j)\le d_{\rm
gr}(\xi_i)$ for all $j\ne i$. If $d_{\rm gr}(\xi_j)<d_{\rm
gr}(\xi_i)$ for all $j\ne i$, then $\xi_i\in\sum^{i-1}_{j=1}A\xi_j$,
which contradicts the assumption. If $d_{\rm gr}(\xi_i)=d_{\rm
gr}(\xi_j)$ for some $j\ne i$, then since $h_j\ne 0$ we have $h_j\in
A_0-\{ 0\}=K^*$. Putting $i'=\max\{i,~j~|~f_j\in K^*\}$, we then
have $\xi_{i'}\in\sum^{i'-1}_{j=1}A\xi_j$, which again contradicts
the assumption. Hence, under the assumption we conclude that $ U$ is 
a minimal homogeneous generating set of $N$.\par

(ii) In view of (i), it is sufficient to show that $\OV U$ is a
homogeneous generating set of $N$. Indeed, if $\xi_i\in  U-\OV
 U$, then $\xi_i\in\sum^{i-1}_{j=1}A\xi_j$. By checking
$\xi_{i-1}$ and so on, it follows that $\xi_i\in\sum_{\xi_k\in\OV 
U}A\xi_k$, as desired.\par\QED{\parindent=0pt\v5

{\bf 3.7. Corollary} Let $ U =\{ \xi_1,\ldots ,\xi_m\}$ be a minimal 
homogeneous generating set of a graded submodule $N$ of $L$, where 
$d_{\rm gr}(\xi_1)\le d_{\rm gr}(\xi_2)\le \cdots\le d_{\rm 
gr}(\xi_m)$, and let $\xi\in L-N$ be a homogeneous element with $ 
d_{\rm gr}(\xi_m)\le d_{\rm gr}(\xi )$. Then $\widehat{ U}= 
 U\cup\{\xi\}$ is a minimal homogeneous generating set of the 
graded submodule $\widehat{N}=N+A\xi$.\QED\v5

{\bf 3.8. Theorem} (Compare with ([KR2], Theorem 4.6.3)) Let $ U =\{ 
\xi_1,\ldots ,\xi_m\}\subset L$ be a finite set of nonzero 
homogeneous elements of $L$ with $d_{\rm gr}(\xi_1)\le d_{\rm 
gr}(\xi_2)\le\cdots \le d_{\rm gr}(\xi_m)$. Then the following 
algorithm  returns a minimal homogeneous generating set $ 
U_{\min}=\{\xi_{j_1},\ldots ,\xi_{j_r}\}\subset  U$ for the graded 
submodule $N=\sum^m_{i=1}A\xi_i$; and meanwhile it returns a 
homogeneous left Gr\"obner basis $\G =\{ g_1,...,g_t\}$ for $N$ such 
that $d_{\rm gr}(g_1)\le d_{\rm gr}(g_2)\le\cdots d_{\rm 
gr}(g_t)$.\vskip 6pt

\underline{\bf Algorithm 3
~~~~~~~~~~~~~~~~~~~~~~~~~~~~~~~~~~~~~~~~~~~~~~~~~~~~~~~~~~~~~~~~~~~~~~~~~~~~~~~~~~~~~~~~~~~~~~~~~}\par

$\begin{array}{l} \textsc{INPUT}: ~ U = \{ 
\xi_1,...,\xi_m\}~\hbox{with}~d_{\rm gr}(\xi_1)\le
d_{\rm gr}(\xi_2)\le\cdots \le d_{\rm gr}(\xi_m)\\
\textsc{OUTPUT}: ~ U_{\min}=\{\xi_{j_1},\ldots ,\xi_{j_r}\}\subset
 U
~\hbox{a minimal homogeneous generating set}\\
\hskip 6.2truecm\hbox{of }~N;\\
\hskip 2.5truecm \G =\{ g_1,...,g_t\}~\hbox{a homogeneous left 
Gr\"obner basis
of}~N\\
\textsc{INITIALIZATION}: ~\S :=\emptyset ,~W := U,~\G :=\emptyset 
,~t':=0,~ U_{\min}:=\emptyset\end{array}$\par 

$\begin{array}{l}
\textsc{LOOP}\\
n :=\min\{d_{\rm gr}(\xi_i),~d_{\rm 
gr}(S_{\ell}(g_i,g_j))~|~\xi_i\in
W,~S_{\ell}(g_i,g_j)\in\S\}\\
\S_n:=\{ S_{\ell}(g_i,g_j)\in\S~|~d_{\rm gr}(S_{\ell}(g_i,g_j))=n\} 
,~
W_n:=\{\xi_j\in W~|~d_{\rm gr}(\xi_j)=n\}\\
\S :=\S -\S_n,~W:=W-W_n\\
~~~~~~\textsc{WHILE}~\S_n\ne\emptyset~\textsc{DO}\\
~~~~~~~~~~~~~\hbox{Choose any}~S_{\ell}(g_i,g_j)\in\S_{n}\\
~~~~~~~~~~~~~\S_n :=\S_n -\{ S_{\ell}(g_i,g_j)\}\\
~~~~~~~~~~~~~\OV{S_{\ell}(g_i,g_j)}^{\G}=\eta\\
~~~~~~~~~~~~~\textsc{IF}~\eta\ne 0~\hbox{with}~\LM (\eta )=
a^{\rho}e_k~\textsc{THEN}\\
~~~~~~~~~~~~~~~~~~t':=t'+1,~g_{t'}:=\eta\\
~~~~~~~~~~~~~~~~~~\S :=\S\cup\{ S_{\ell}(g_i,g_{t'})~|~g_i\in\G ,~1\le i< t',~\LM (g_i )=a^{\tau}e_k\}\\
~~~~~~~~~~~~~~~~~~\G :=\G\cup\{ g_{t'}\}\\
~~~~~~~~~~~~~\textsc{END}\\
~~~~~~\textsc{END}\\
~~~~~~\textsc{WHILE}~W_n\ne\emptyset~\hbox{DO}\\
~~~~~~~~~~~~~\hbox{Choose any}~\xi_j\in W_{n}\\
~~~~~~~~~~~~~W_n :=W_n -\{ \xi_j\}\\
~~~~~~~~~~~~~\OV{\xi_j}^{\G}=\eta\\
~~~~~~~~~~~~~\textsc{IF}~\eta\ne 0~\hbox{with}~\LM (\eta )=
a^{\rho}e_k~\textsc{THEN}\\
~~~~~~~~~~~~~~~~~~U_{\min}:= U_{\min}\cup \{\xi_j\}\\
~~~~~~~~~~~~~~~~~~t':=t'+1,~g_{t'}:=\eta\\
~~~~~~~~~~~~~~~~~~\S :=\S\cup\{ S_{\ell}(g_i,g_{t'})~|~g_i\in\G,~1\le i<t',~\LM (g_i )=a^{\tau}e_k\}\\
~~~~~~~~~~~~~~~~~~\G :=\G\cup\{ g_{t'}\}\\
~~~~~~~~~~~~~\textsc{END}\\
\end{array}$

$\begin{array}{l}
~~~~~~\textsc{END}\\
\textsc{UNTIL}~\S=\emptyset\\
\textsc{END}\\
\end{array}$\par
\underline{~~~~~~~~~~~~~~~~~~~~~~~~~~~~~~~~~~~~~~~~~~~~~~~~~~~~~~~~~~~~~~~~~~~~~~~~~~~~~~~~~~~~~~~~~~~~~~~~~~~~~~~~~~~~~~~~~~~~~~~~~~~~~~~~~~~~}
\vskip 6pt

{\bf Proof} Since this algorithm is clearly a variant of {\bf
Algorithm 1} and {\bf Algorithm 2} with a minimization procedure, it
terminates and the eventually obtained $\G$ is a homogeneous left
Gr\"obner basis for $N$ in which the degrees of elements are ordered
non-decreasingly. It remains to prove that the eventually obtained $ 
U_{\min}$ is a minimal homogeneous generating set of $N$.}\par

As in the proof of Proposition 3.5, let us first bear in mind that
for each $n$, in both the WHILE loops every  newly appended
$S_{\ell}(g_i,g_{t'})$ has $d_{\rm gr}(S_{\ell}(g_i,g_{t'}))>n$.
Moreover, for  convenience, let us write $\G (n)$ for the $\G$
obtained after $\S_n$ is exhausted in the first WHILE loop, and
write $ U_{\min}[n]$, $\G [n]$ respectively  for the $ U_{\min}$, 
$\G$ obtained after $W_n$ is exhausted in the second WHILE loop.  
Since the algorithm starts with ${\cal O}=\emptyset$ and  
$\G=\emptyset$, if, for a fixed $n$, we check carefully how the 
elements of $ U_{\min}$ are chosen during executing the second WHILE 
loop, and how the new elements are appended to $\G$ after each pass 
through the first or the second WHILE loop, then it follows from 
Proposition 3.3 and Corollary 3.4 that after $W_n$ is exhausted, the 
obtained $ U_{\min}[n]$ and $\G [n]$ generate the same module, 
denoted $N(n)$, such that  $\G [n]$ is an $n$-truncated left 
Gr\"obner basis of $N(n)$.  We now use induction to show that the 
eventually obtained  $ U_{\min}$ is a minimal homogeneous generating 
set for $N$. If $ U_{\min}=\emptyset$, then it is a minimal 
generating set of the zero module.  To proceed,  we assume that $ 
U_{\min}[n]$ is a minimal homogeneous generating set for $N(n)$ 
after $W_n$ is exhausted in the second WHILE loop. Suppose that 
$n_1$ is the first number after $n$ such that 
$\S_{n_1}\ne\emptyset$. We complete the induction proof below by 
showing  that $ U_{\min}[n_1]$ is a minimal homogeneous generating 
set of $N(n_1)$.\par

If in a certain pass through the first WHILE loop,
$\OV{S_{\ell}(g_i,g_j)}^{\G}=\eta\ne 0$ for some
$S_{\ell}(g_i,g_j)\in\S_{n_1}$, then we note that $\eta\in N(n)$. It
follows that after $\S_{n_1}$ is exhausted in the first WHILE loop,
$\G (n_1)$ generates $N(n)$ and $\G (n_1)$ is an $n_1$-truncated
left Gr\"obner basis of $N(n)$. Next, assume that $W_{n_1}=\{
\xi_{j_1},\ldots ,\xi_{j_s}\}\ne\emptyset$ and that the elements of
$W_{n_1}$ are processed in the given order during executing the
second WHILE loop. Since $\G (n_1)$ is an $n_1$-truncated left
Gr\"obner basis of $N(n)$, if $\xi_{j_1}\in W_{n_1}$ is such that
$\OV {\xi_{j_1}}^{\G (n_1)}=\eta_1\ne 0$, then $\xi_{j_1}, \eta_1\in
L -N(n)$. By Corollary 3.4, we conclude  that $\G
(n_1)\cup\{\eta_1\}$ is an $n_1$-truncated Gr\"obner basis for the
module $N(n)+A\eta_1$; and by Corollary 3.7, we conclude that $ 
U_{\min}[n]\cup\{\xi_{j_1}\}$ is a minimal homogeneous generating 
set of $N(n)+A\eta_1$. Repeating this procedure, if $\xi_{j_2}\in 
W_{n_1}$ is such that $\OV{f_{j_2}}^{\G 
(n_1)\cup\{\eta_1\}}=\eta_2\ne 0$, then $\xi_{j_2}, \eta_2\in L 
-(N(n)+A\eta_1)$.  By Corollary 3.4, we conclude  that $\G 
(n_1)\cup\{\eta_1,\eta_2\}$ is an $n_1$-truncated left Gr\"obner 
basis for the module  $N(n)+A\eta_1+A\eta_2$; and by Corollary 3.7, 
we conclude that $ U_{\min}[n]\cup\{\xi_{j_1},\xi_{j_2}\}$ is a 
minimal homogeneous generating set of $N(n)+A\eta_1+A\eta_2$. 
Continuing this procedure until $W_{n_1}$ is exhausted we assert 
that the returned $\G [n_1]=\G$ and $ U_{\min}[n_1]= U_{\min}$ 
generate the same module $N(n_1)$ and $\G [n_1]$ is an 
$n_1$-truncated left Gr\"obner basis of $N(n_1)$ and $ 
U_{\min}[n_1]$ is a minimal homogeneous generating set of $N(n_1)$, 
as desired. As all elements of $ U$ are eventually processed by the 
second WHILE loop, we conclude that the finally obtained $\G$ and $ 
U_{\rm min}$ have the properties: $\G$ generates the module $N$, 
$\G$ is an $n_0$-truncated left Gr\"obner basis of $N$, and $ 
U_{\min}$ is a minimal homogeneous generating set of $N$.\par\QED 
{\parindent=0pt \v5

{\bf Remark} If we are only interested in getting a minimal 
homogeneous generating set for the submodule $N$, then {\bf 
Algorithm 3} can indeed be speed up. More precisely, with  $$d_{\rm 
gr}(\xi_1)\le d_{\rm gr}(\xi_2)\le\cdots \le d_{\rm gr}(\xi_m)= 
n_0,$$ it follows from the proof above that if we  stop executing 
the algorithm after $S_{n_0}$ and $W_{n_0}$ are exhausted, then the 
resulted $ U_{\min}[n_0]$ is already the desired minimal homogeneous 
generating set for $N$, while $\G [n_0]$ is an $n_0$-truncated left 
Gr\"obner basis of $N$. \v5

{\bf 3.9. Corollary}  Let $ U =\{ \xi_1,\ldots ,\xi_m\}\subset L$ be 
a finite set of nonzero homogeneous elements of $L$ with $d_{\rm 
gr}(\xi_1)=d_{\rm gr}(\xi_2)=\cdots = d_{\rm gr}(\xi_m)=n_0$
\par

(i) If $ U$ satisfies $\LM (\xi_i)\ne\LM (\xi_j)$ for all $i\ne j$, 
then $ U$ is a minimal homogeneous generating set of the graded 
submodule $N=\sum_{i=1}^mA\xi_i$ of $L$, and meanwhile $ U$ is an 
$n_0$-truncated left Gr\"obner basis for $N$.\par

(ii) If $ U$ is a minimal left Gr\"obner basis of the graded
submodule $N=\sum^m_{i=1}A\xi_i$ (i.e., $ U$ is a left Grobner basis 
of $N$ satisfying $\LM (\xi_i){\not |}~\LM (\xi_j)$ for all $i\ne 
j$), then  $ U$ is a minimal homogeneous generating set of 
$N$.\vskip 6pt

{\bf Proof} By the assumption, it follows from the second WHILE loop
of {\bf Algorithm 3} that $ U_{\min}= U$. }\v5

\section*{4. Computation of Minimal Homogeneous Generating Sets of Graded Quotient Modules $M=L/N$}
In this section, $A=K[a_1,\ldots ,a_n]$ denotes a weighted
$\NZ$-graded solvable polynomial algebra with respect to a given
positive-degree function $d(~)$ on $A$, $(\B,\prec )$ denotes a
fixed admissible system of $A$, and $L_0=\oplus_{i=1}^sAe_i$ denotes
a graded free left $A$-module such that $d_{\rm gr}(e_i)=b_i$, $1\le
i\le s$, i.e., $L_0=\oplus_{q\in\NZ}L_{0q}$ with $L_{0q}=K$-span$\{
a^{\alpha}e_i\in\BE~|~d(a^{\alpha})+b_i=q\}$. Let $N$ be a graded
submodule of $L_0$ and $M=L_0/N$. By mimicking ([KR], Proposition
4.7.24)), our aim of this section is to present an algorithm for
computing a minimal homogeneous generating set of the graded
$A$-module $M$. All conventions and notations used before are
maintained.\v5

Since $A$ is Noetherian, $N$ is a finitely generated graded
submodule of $L_0$. Let $N=\sum^m_{j=1}A\xi_j$ be generated by the
set of nonzero homogeneous elements $ U =\{\xi_1,\ldots ,\xi_m\}$, 
where $\xi_{\ell}=\sum_{k=1}^sf_{k\ell}e_k$ with $f_{k\ell}\in A$, 
$1\le \ell\le m$. Then, every nonzero $f_{k\ell}$ is a homogeneous 
element of $A$ such that $d_{\rm gr}(\xi_{\ell})=d_{\rm 
gr}(f_{k\ell}e_k)=d_{\rm gr}(f_{k\ell})+b_k$, where $b_k=d_{\rm 
gr}(e_k)$, $1\le k\le s$, $1\le \ell\le m$. {\parindent=0pt\v5

{\bf 4.1. Lemma} With  every $\xi_{\ell}=\sum_{i=1}^sf_{i\ell}e_i$
as fixed above, $1\le \ell\le m$, if the $i$-th coefficient $f_{ij}$
of some $\xi_j$ is a nonzero constant, say $f_{ij}=1$ without loss
of generality, then for each $\ell =1,\ldots ,j-1,j+1,\ldots ,m$,
the element $\xi_{\ell}'=\xi_{\ell}-f_{i\ell}\xi_j$ does not involve
$e_i$. Putting $ U '=\{ \xi_1',\ldots ,\xi'_{j-1},\xi_{j+1}',\ldots 
,\xi'_m\}$, there is a graded $A$-module isomorphism 
$M'=L_0'/N'\cong L_0/N=M$, where $L_0'=\oplus_{k\ne i}Ae_k$ and 
$N'=\sum_{\xi_{\ell}'\in U '}A\xi_{\ell}'$. \vskip 6pt

{\bf Proof} Since $f_{ij}=1$ by the assumption, we see that every
$\xi_{\ell}'=\sum_{k\ne i}(f_{k\ell}-f_{i\ell}f_{kj})e_k$ does not
involve $e_i$. Let $ U '=\{ \xi_1',\ldots
,\xi'_{j-1},\xi_{j+1}',\ldots ,\xi'_m\}$ and
$N'=\sum_{\xi_{\ell}'\in U '}A\xi_{\ell}'$. Then $N'\subset
L_0'=\oplus_{k\ne i}Ae_k$. Again since $f_{ij}=1$, we have $d_{\rm
gr}(\xi_j)=d_{\rm gr}(e_i)=b_i$. It follows from the property
($\mathbb{P}$4) formulated in Section 3 that
$$\begin{array}{rcl} d_{\rm gr}(f_{i\ell}f_{kj}e_k)&=&d_{\rm gr}(f_{i\ell})+d_{\rm gr}(f_{kj}e_k)\\
&=&d_{\rm gr}(f_{i\ell})+d_{\rm gr}(\xi_j)\\
&=&d_{\rm gr}(f_{i\ell})+b_i\\
&=&d_{\rm gr}(f_{i\ell}e_i)\\
&=&d_{\rm gr}(\xi_{\ell})\\
&=&d_{\rm gr}(f_{k\ell}e_k).\end{array}$$ Noticing that $d_{\rm
gr}(f_{i\ell}\xi_j)=d_{\rm gr}(f_{i\ell})+d_{\rm gr}(\xi_j)$, this
shows that in the representation of $\xi_{\ell}'$ every nonzero term
$(f_{k\ell}-f_{i\ell}f_{kj})e_k$ is a homogeneous element of degree
$d_{\rm gr}(\xi_{\ell})=d_{\rm gr}(f_{i\ell}\xi_j)$, thereby
$M'=L_0'/N'$ is a graded $A$-module. Note that $N=N'+A\xi_j$ and
that  $\xi_j =e_i+\sum_{k\ne i}f_{kj}e_k$. Without making confusion,
if we use the same notation $\OV e_k$ to denote the coset
represented by $e_k$ in $M'$ and $M$ respectively,  it is now clear
that the desired graded $A$-module isomorphism
$M'\mapright{\varphi}{}M$ is naturally defined by $\varphi (\OV
e_k)=\OV e_k$, $k=1,\ldots ,i-1,i+1,\ldots ,s$. \QED}\v5

Let $M=L_0/N$ be as fixed above with $N$ generated by the set of 
nonzero homogeneous elements $ U =\{ \xi_1,\ldots ,\xi_m\}$. Then 
since $A$ is $\NZ$-graded with $A_0=K$, from the literature (e.g. 
([Eis], Chapter 19), ([Kr1], Chapter 3), [Li3]) it is well known 
that the homogeneous generating set $\OV E=\{ \OV e_1,\ldots ,\OV 
e_s\}$ of $M$ is a minimal homogeneous generating set if and only if 
$\xi_{\ell}=\sum_{k=1}^sf_{k\ell}e_k$ implies $d_{\rm 
gr}(f_{k\ell})>0$ whenever $f_{k\ell}\ne 0$, $1\le\ell\le m$. 
{\parindent=0pt\v5

{\bf 4.2. Proposition} (Compare with ([KR2], Proposition 4.7.24)) 
With notation as fixed above, the algorithm presented below returns 
a subset $\{ e_{i_1},\ldots ,e_{i_{s'}}\}\subset\{ e_1,\ldots 
,e_s\}$ and a subset $V=\{ v_1,\ldots ,v_t\}\subset N\cap L_0'$ such 
that $M\cong L_0'/N'$ as graded $A$-modules, where 
$L_0'=\oplus_{q=1}^{s'}Ae_{i_q}$ with $s'\le s$ and 
$N'=\sum^{t}_{k=1}Av_k$, and such that $\{\OV e_{i_1},\ldots ,\OV 
e_{i_{s'}}\}$ is a minimal homogeneous generating set of $M$. \vskip 
6pt

\underline{\bf Algorithm 4
~~~~~~~~~~~~~~~~~~~~~~~~~~~~~~~~~~~~~~~~~~~~~~~~~~~~~~~~~~~~~~~~~~~~~~~~~~~~~~~~~~~~~~~~~~~~~~~~~}\par

$\begin{array}{l} \textsc{INPUT}: E=\{ e_1,\ldots ,e_s\};~~ U =
\{ \xi_1,...,\xi_m\}\\
~~~~~~~~~~~~~\hbox{where}~\xi_{\ell}=\sum_{k=1}^sf_{k\ell}e_k~
\hbox{with homogeneous}~f_{k\ell}\in A, ~1\le \ell\le m.\\
\textsc{OUTPUT}: ~E' =\{ e_{i_1},\ldots ,e_{i_{s'}}\}; ~~V=\{
v_1,\ldots ,v_t\}\subset N\cap
L_0',~\\
~~~~~~~~~~~~~\hbox{such that}~v_{j}=\sum^{s'}_{q=1}h_{qj}e_{i_q}\in
L_0'=\oplus_{q=1}^{s'}Ae_{i_q}~
\hbox{with}~h_{qj}\not\in K^*\\
~~~~~~~~~~~~~\hbox{whenever}~h_{qj}\ne 0,~1\le j\le t\\
\textsc{INITIALIZATION}:    ~t :=m;~ V := U ;  ~s':=s; ~E' :=E;
\end{array}$\par

$\begin{array}{l}
\textsc{BEGIN}\\
~~~~\textsc{WHILE}~\hbox{there is a}~v_j=\sum^{s'}_{k=1}f_{kj}e_k\in V~\hbox{satisfying}\\
~~~~~~~~~~~~~~~~~f_{kj}\not\in K^*~\hbox{for}~k<i~\hbox{and}~f_{ij}\in K^*~\hbox{DO}\\
~~~~~~~~~~~~~~~~~\hbox{for}~T=\{1,\ldots ,j-1,j+1,\ldots
,t\}~\hbox{compute}\\
~~~~~~~~~~~~~~~~~v_{\ell}'
=v_{\ell}-\frac{1}{f_{ij}}f_{i\ell}v_j,~\ell\in T,~r=\#\{\ell~|~\ell\in T,~v_{\ell}' =0\}\\
~~~~~~~~~~t := t-r-1\\
~~~~~~~~~V :=\{ v_{\ell}=v_{\ell}'~|~\ell\in T,~v_{\ell}'\ne 0\}\\
~~~~~~~~~~~~=\{ v_1,\ldots ,v_t\}~(\hbox{after reordered})\\
~~~~~~~~~s' :=s'-1\\
~~~~~~~~~E':=E'-\{ e_i\}=\{ e_1,\ldots ,e_{s'}\}~(\hbox{after reordered})\\
~~~~\textsc{END}\\
\textsc{END}\end{array}$\par
\underline{~~~~~~~~~~~~~~~~~~~~~~~~~~~~~~~~~~~~~~~~~~~~~~~~~~~~~~~~~~~~~~~~~~~~~~~~~~~~~~~~~~~~~~~~~~~~~~~~~~~~~~~~~~~~~~~~~~~~~~~~~~~~~~~~~~}
\vskip 6pt

{\bf Proof} It is clear that the algorithm is finite. The
correctness of the algorithm follows immediately from Lemma 4.1 and
the remark we made before the proposition. }\v5

\section*{5. Computation of Minimal Finite Graded Free Resolutions}
Let $A=K[a_1,\ldots ,a_n]$ be a weighted  $\NZ$-graded solvable 
polynomial algebra with respect to a given positive-degree function 
$d(~)$ on $A$, and $(\B,\prec )$ a fixed admissible system of $A$. 
Then since $A$ is Noetherian and $A_0=K$, it is theoretically well 
known that up to a graded isomorphism of chain complexes in the 
category of graded $A$-modules, every finitely generated graded 
$A$-module $M$ has a unique minimal graded free resolution (cf. 
([Eis], Chapter 19), ([Kr1], Chapter 3), [Li3]). Combining the 
results of Section 3 and Section 4, in this section we present 
algorithmic procedures for constructing  such a minimal graded free 
resolution over $A$. \par

All notions, notations and conventions introduced before are 
maintained.\v5

In what follows, $L=\oplus_{i=1}^sAe_i$ denotes a graded free left
$A$-module such that $d_{\rm gr}(e_i)=b_i$, $1\le i\le s$, i.e.,
$L=\oplus_{q\in\NZ}L_{q}$ with $L_{q}=K$-span$\{
a^{\alpha}e_i\in\BE~|~d(a^{\alpha})+b_i=q\}$, and $\prec_{e}$ 
denotes a  left monomial ordering on the $K$-basis $\BE$ of $L$. 
Moreover,  as in Section 2 we write $S_{\ell}(\xi_i,\xi_j)$ for the 
left S-polynomial of two elements $\xi_i$, $\xi_j\in L$.
\par

Let $N=\sum^m_{i=1}A\xi_i$ be a graded submodule of $L$ generated by
the set of nonzero homogeneous elements $ U =\{\xi_1,\ldots
,\xi_m\}$. We first demonstrate how to calculate a generating set of
the syzygy module Syz$( U )$  by means of a left Gr\"obner basis of 
$N$. To this end, let $\G=\{ g_1,\ldots ,g_t\}$ be a left Gr\"obner 
basis of $N$ with respect to $\prec_{e}$, then every nonzero left 
S-polynomial $S_{\ell}(g_i,g_j)$ has a left Gr\"obner representation 
$S_{\ell}(g_i,g_j)=\sum_{i=1}^tf_ig_i$ with 
$\LM(f_ig_i)\preceq_{\mathbbm{e}}\LM (S_{\ell}(g_i,g_j))$ whenever 
$f_i\ne 0$ (note that such a representation is obtained by using the 
division by $\G$ during executing the WHILE loop in {\bf Algorithm 
1} or {\bf Algorithm 3}). Considering the syzygy module Szy$(\G )$ 
of $\G$ in the free $A$-module $L_1=\oplus_{i=1}^tA\varepsilon_i$, 
if we put
$$s_{ij}=f_1\varepsilon_1+\cdots +\left (f_i-\frac{a^{\gamma-\alpha (i)}}{\LC (a^{\gamma -\alpha (i)}\xi_i)}\right )\varepsilon_i+
\cdots +\left (f_j+\frac{a^{\gamma -\alpha (j)}}{\LC (a^{\gamma 
-\alpha (j)}\xi_j)}\right )\varepsilon_j+\cdots +f_t\varepsilon_t,$$ 
${\cal S} =\{ s_{ij}~|~1\le i<j\le t\}$, then it can be shown, 
actually as in the commutative case (cf. [AL], Theorem 3.7.3), that 
${\cal S}$ generates Szy$(\G )$. However, by employing an analogue 
of the Schreyer ordering $\prec_{s\hbox{-}\varepsilon}$ on the 
$K$-basis $\B (\varepsilon )=\{ 
a^{\alpha}\varepsilon_i~|~a^{\alpha}\in\B ,~1\le i\le m\}$ of $L_1$ 
induced by $\G$ with respect to $\prec_{e}$, which is defined 
subject to the rule: for 
$a^{\alpha}\varepsilon_i,a^{\beta}\varepsilon_j\in\B (\varepsilon 
)$,
$$a^{\alpha}\varepsilon_i\prec_{s\hbox{-}\varepsilon} a^{\beta}\varepsilon_j
\Leftrightarrow\left\{\begin{array}{l}
\LM (a^{\alpha}g_i)\prec_{e}\LM (a^{\beta}g_j),\\
\hbox{or}\\
\LM (a^{\alpha}g_i)=\LM (a^{\beta}g_j)~\hbox{and}~i<j,
\end{array}\right.$$
there is indeed a much stronger result, namely the noncommutative
analogue of Schreyer's Theorem [Sch] (cf. Theorem 3.7.13 in [AL] for
free modules over commutative polynomial algebras; Theorem 4.8 in
[Lev] for free modules over solvable polynomial
algebras):{\parindent=0pt \v5

{\bf 5.1. Theorem}  With respect to the left monomial ordering 
$\prec_{s\hbox{-}\varepsilon}$ on $\B (\varepsilon )$ as defined 
above, the following statements hold.\par

(i) Let $s_{ij}$ be determined by $S_{\ell}(g_i,g_j)$, where $i<j$,
$\LM (g_i)=a^{\alpha (i)}e_s$ with $\alpha (i)=(\alpha_{i_1},\ldots
,\alpha_{i_n})$, and $\LM (g_j)=a^{\alpha (j)}e_s$ with $\alpha
(j)=(\alpha_{j_1},\ldots ,\alpha_{j_n})$. Then $\LM
(s_{ij})=a^{\gamma-\alpha (j)}\varepsilon_j$, where $\gamma
=(\gamma_1,\ldots ,\gamma_n)$ with each
$\gamma_k=\max\{\alpha_{i_k},\alpha_{j_k}\}$. \par

(ii) ${\cal S}$ is a left Gr\"obner basis of Syz$(\G )$, thereby
${\cal S}$ generates Syz$(\G )$.\par\QED} \v5

To go further, again let $\G=\{ g_1,\ldots ,g_t\}$ be  the left
Gr\"obner basis of $N$ produced by running {\bf Algorithm 1} (or
{\bf Algorithm 3}) with the initial input data $ U=\{\xi_1,\ldots 
,\xi_m\}$. Using the usual matrix notation for convenience, we have
$$\left (\begin{array}{l} \xi_1\\ \vdots\\ \xi_m\end{array}\right )=U_{m\times t}\left (\begin{array}{l} g_1
\\ \vdots\\ g_t\end{array}\right ) ,\quad
\left (\begin{array}{l} g_1\\ \vdots\\ g_t\end{array}\right
)=V_{t\times m}\left (\begin{array}{l} \xi_1
\\ \vdots\\ \xi_m\end{array}\right ) ,$$
where the $m\times t$ matrix $U_{m\times t}$ (with entries in $A$)
is obtained by the division by $\G$, and the $t\times m$ matrix
$V_{t\times m}$ (with entries in $A$) is obtained by keeping track
of the reductions during executing the WHILE loop of {\bf Algorithm
1}.  By Theorem 5.1, we may write Syz$(\G )=\sum^r_{i=1}A\mathcal
{S}_i$ with $\mathcal{S}_1,\ldots ,\mathcal{S}_r\in
L_1=\oplus_{i=1}^tA\varepsilon_i$;  and if
$\mathcal{S}_i=\sum^t_{j=1}f_{ij}\varepsilon_j$, then we  write
$\mathcal{S}_i$ as a $1\times t$ row matrix, i.e.,
$\mathcal{S}_i=(f_{i_1}~\ldots~ f_{it})$,  whenever matrix notation
is convenient in the according discussion. At this point, we note
also that all the $\mathcal{S}_i$ may be written down one by one
during executing the WHILE loop of {\bf Algorithm 1} (or the first
WHILE loop in {\bf Algorithm 3}) successively. Furthermore, we write
$D_{(1)},\ldots ,D_{(m)}$ for the rows of the matrix $D_{m\times
m}=U_{m\times t}V_{t\times m}-E_{m\times m}$ where $E_{m\times m}$
is the $m\times m$ identity matrix. The following proposition is a
noncommutative analogue of ([AL],  Theorem 3.7.6). {\parindent=0pt 
\v5

{\bf Proposition 5.2.}  With notation fixed above, the syzygy module
Syz$( U )$ of $ U =\{ \xi_1,\ldots ,\xi_m\}$ is generated by
$$\{ \mathcal{S}_1V_{t\times m},\ldots ,\mathcal{S}_{r}V_{t\times m},D_{(1)},\ldots ,D_{(m)}\} ,$$
where each $1\times m$ row matrix represents an element of the free
$A$-module $\oplus_{i=1}^mA\omega_i$. \vskip 6pt

{\bf Proof} Since $$0=\mathcal{S}_i\left (\begin{array}{l} g_1\\
\vdots\\ g_t\end{array}\right )=(f_{i_1}~\ldots ~f_{it})\left
(\begin{array}{l} g_1\\ \vdots\\ g_t\end{array}\right
)=(f_{i_1}~\ldots ~f_{it})V_{t\times m}\left (\begin{array}{l}
\xi_1\\ \vdots\\ \xi_m\end{array}\right ) ,$$ we have
$\mathcal{S}_iV_{t\times m}\in$ Syz$( U )$, $1\le i\le r$. Moreover, 
since
$$\begin{array}{rcl} D_{m\times m}\left (\begin{array}{l}
\xi_1\\ \vdots\\ \xi_m\end{array}\right )&=&(U_{m\times t}V_{t\times
m}-E_{m\times m})\left (\begin{array}{l} \xi_1\\ \vdots\\
\xi_m\end{array}\right )\\
\\
&=&U_{m\times t}V_{t\times m}\left (\begin{array}{l} \xi_1\\
\vdots\\ \xi_m\end{array}\right )-\left (\begin{array}{l} \xi_1\\
\vdots\\ \xi_m\end{array}\right )\\
\\
&=&U_{m\times t}\left (\begin{array}{l} g_1\\ \vdots\\
g_t\end{array}\right )-\left (\begin{array}{l} \xi_1\\ \vdots\\
\xi_m\end{array}\right )=\left (\begin{array}{l} \xi_1\\ \vdots\\
\xi_m\end{array}\right )-\left (\begin{array}{l} \xi_1\\ \vdots\\
\xi_m\end{array}\right )=0,\end{array}$$ we have $D_{(1)},\ldots
,D_{(r)}\in$ Syz$( U )$.}\par

On the other hand, if $H =(h_1~\ldots~h_m)$ represents the element
$\sum^m_{i=1}h_i\omega_i\in \oplus_{i=1}^mA\omega_i$ such that
$H\left (\begin{array}{l} \xi_1\\ \vdots\\ \xi_m\end{array}\right
)=0$, then $0=H U_{m\times t}\left (\begin{array}{l} g_1\\ \vdots\\
g_t\end{array}\right )$. This means $H U_{m\times t}\in$ Syz$(\G )$.
Hence, $H U_{m\times t}=\sum^r_{i=1}f_i\mathcal{S}_i$ with $f_i\in
A$, and it follows that $H U_{m\times t}V_{t\times
m}=\sum^r_{i=1}f_i\mathcal{S}_iV_{t\times m}$. Therefore,
$$\begin{array}{rcl} H&=&H +H U_{m\times t}V_{t\times m}-H U_{m\times t}V_{t\times m}\\
&=&H (E_m-U_{m\times t}V_{t\times
m})+\sum^r_{i=1}f_i\mathcal{S}_iV_{t\times m}\\
&=&-H D_{m\times m}+\sum^r_{i=1}f_i(\mathcal{S}_iV_{t\times
m}).\end{array}$$ This shows that every element of Syz$( U )$ is
generated by $\{ \mathcal{S}_1V_{t\times m},\ldots
,\mathcal{S}_{r}V_{t\times m},D_{(1)},\ldots ,D_{(m)}\} ,$ as
desired.\QED\v5

Next,  we recall the noncommutative version of Hilbert's syzygy
theorem for solvable polynomial algebras. For a constructive proof
of Hilbert's syzygy theorem by means of Gr\"obner bases respectively
in  the commutative case and the noncommutative case, we refer to
(Corollary 15.11 in  [Eis]) and (Section 4.4 in [Lev]).
{\parindent=0pt\v5

{\bf 5.3. Theorem} Let $A=K[a_1,\ldots ,a_n]$ be a solvable
polynomial algebra with admissible system $(\B ,\prec )$. Then every
finitely generated left $A$-module $M$ has a free resolution
$$0~\mapright{}{}~L_s~\mapright{}{}~L_{s-1}~\mapright{}{}\cdots ~
\mapright{}{}~L_0~\mapright{}{}~M~\mapright{}{}~ 0$$ where each
$L_i$ is a free $A$-module of finite rank and $s\le n$. It follows
that $M$ has projective dimension p.dim$_AM\le s$, and that $A$ has
global homological dimension gl.dim$A\le n$. \par\QED}\v5

Now, we are able to give the main result of this section. Let
$M=\sum^s_{i=1}Av_i$ be a finitely generated $\NZ$-graded left
$A$-module with the set of homogeneous generators $\{ v_1,\ldots
,v_s\}$ such that $d_{\rm gr}(v_i)=b_i$ for $1\le i\le s$, i.e.,
$M=\oplus_{b\in\NZ}M_b$ with each $M_b$ the degree-$b$ homogeneous
part. Then $M$ is  isomorphic to a quotient module of the
$\NZ$-graded free $A$-module $L_0=\oplus_{i=1}^sAe_i$ which is
equipped with the $\NZ$-gradation  $L_0=\oplus_{q\in\NZ}L_{0q}$ such
that $d_{\rm gr}(e_i)=b_i$, $1\le i\le s$, and $L_{0q}=K$-span$\{
a^{\alpha}e_i\in\BE~|~d(a^{\alpha})+b_i=q\}$, where $d(~)$ is the
given positive-degree function on $A$. Thus we may write $M=L_0/N$,
where $N$ is a graded submodule of $L_0$. \par

Recall that a {\it minimal graded free resolution} of $M$ is an
exact sequence by free $A$-modules and  $A$-module homomorphisms
$${\cal L}_{\bullet}\quad\quad\cdots~\mapright{\varphi_{i+1}}{}~L_i~\mapright{\varphi_i}{}~\cdots ~
\mapright{\varphi_2}{}~L_1~\mapright{\varphi_1}{}~L_0~\mapright{\varphi_0}{}~M~\mapright{}{}~0$$
in which each $L_i$ is an $\NZ$-graded free $A$-module with a finite
homogeneous $A$-basis $E_i=\{ e_{i_1},\ldots ,e_{i_{s_i}}\}$, and
each $\varphi_i$ is a graded $A$-module homomorphism of degree 0
(i.e., $\varphi_i$ sends the degree-$q$ homogeneous part of $L_i$
into the degree-$q$ homogeneous part of $L_{i-1}$), such
that{\parindent=1.2truecm\par

\item{(1)} $\varphi_0(E_0)$ is a minimal homogeneous
generating set of $M$, Ker$\varphi_0=N$, and\par

\item{(2)} for $i\ge 1$,  $\varphi_i(E_i)$ is a minimal
homogeneous generating set of
Ker$\varphi_{i-1}$.\par}{\parindent=0pt\v5

{\bf 5.4. Theorem} With notation as fixed above, suppose that
$N=\sum^m_{i=1}A\xi_i$ with the set of nonzero homogeneous
generators $ U =\{\xi_1,\ldots ,\xi_m\}$. Then the graded $A$-module 
$M=L_0/N$ has a minimal graded free resolution of length $d\le n$:
$${\cal L}_{\bullet}\quad\quad 0~\mapright{}{}~L_d~\mapright{\varphi_{q}}{}~
\cdots~\mapright{\varphi_2}{}~L_1~\mapright{\varphi_1}{}~L_0~\mapright{\varphi_0}{}~M~\mapright{}{}~0$$
which can be constructed by implementing the following
procedures:}\par

{\bf Procedure 1.}  Run {\bf Algorithm 4} with the initial input
data $E=\{ e_1,\ldots ,e_s\}$ and $ U =\{ \xi_1,\ldots ,\xi_m\}$ to 
compute a subset $E'=\{ e_{i_1},\ldots ,e_{i_{s'}}\}\subset\{ 
e_1,\ldots ,e_s\}$ and a subset $V=\{ v_1,\ldots ,v_t\}\subset N\cap 
L_0'$ such that $M\cong L_0'/N'$ as graded $A$-modules, where 
$L_0'=\oplus_{q=1}^{s'}Ae_{i_q}$ with $s'\le s$ and 
$N'=\sum^{t}_{k=1}Av_k$, and such that $\{\OV e_{i_1},\ldots ,\OV 
e_{i_{s'}}\}$ is a minimal homogeneous generating set of $M$.\par

For convenience, after accomplishing Procedure 1 we may assume that
$E=E'$, $ U =V$ and $N=N'$. Accordingly we have the short exact
sequence $0~\mapright{}{}~ N~\mapright{}{}~
L_0~\mapright{\varphi_0}{}~M~\mapright{}{}  ~0$ such that $\varphi_0
(E)=\{ \OV e_1,\ldots ,\OV e_s\}$  is a minimal homogeneous
generating set of $M$. \par

{\bf Procedure 2.} Choose a left monomial ordering $\prec_{e}$ on 
the $K$-basis $\BE$ of $L_0$  and run {\bf Algorithm 3} with the 
initial input data $ U =\{ \xi_1,\ldots ,\xi_m\}$ to compute a 
minimal homogeneous generating set $ U_{\min}=\{\xi_{j_1},\ldots 
,\xi_{j_{s_1}}\}$ and a left Gr\"obner basis $\G$ for $N$; at the 
same time, by keeping track of the reductions during executing the 
first WHILE loop and the second WHILE loop respectively, return the 
matrices ${\cal S}$ and $V$ required by Proposition 5.2.
\par

{\bf Procedure 3.} By using the division by the left Gr\"obner basis
$\G$ obtained in Procedure 2, compute the matrix $U$ required by
Proposition 5.2. Use the matrices ${\cal S}$, $V$ obtained in
Procedure 2, the matrix $U$ and Proposition 5.2 to compute a
homogeneous generating set of $N_1=$ Syz$( U_{\min})$ in the
$\NZ$-graded free $A$-module $L_1=\oplus_{i=1}^{s_1}A\varepsilon_i$,
where the gradation of $L_1$ is defined by setting $d_{\rm
gr}(\varepsilon_k)=d_{\rm gr}(\xi_{j_k})$, $1\le k\le s_1$.\par

{\bf Procedure 4.} Construct the exact sequence
$$0~\mapright{}{}~ N_1~\mapright{}{}~L_1~\mapright{\varphi_1}{}~
L_0~\mapright{\varphi_0}{}~M~\mapright{}{}  ~0$$ where
$\varphi_1(\varepsilon_k)=\xi_{j_k}$, $1\le k\le s_1$. If $N_1\ne
0$, then repeat Procedure 2 -- Procedure 4 for $N_1$ and so on.
\par

By Hilbert's syzygy theorem for solvable polynomial algebras 
(Theorem 5.3), $M$ has finite projective dimension p.dim$_AM\le n$. 
Also we know from the literature (e.g. ([Eis], Chapter 19), ([Kr1], 
Chapter 3), [Li3]) that p.dim$_AM=$ the length of the minimal 
resolution of $M$. It follows that the desired minimal graded free 
resolution ${\cal L}_{\bullet}$ for $M$ is then obtained after 
finite times of processing the above procedures.\v5

\v5 \centerline{References}
\parindent=1.6truecm

\re{[AL1]} J. Apel and W. Lassner, An extension of Buchberger's
algorithm and calculations in enveloping fields of Lie algebras.
{\it J. Symbolic Comput}., 6(1988), 361--370.

\item{[AL2]} W. W. Adams and P. Loustaunau, {\it An Introduction to Gr\"obner Bases}.
Graduate Studies in Mathematics, Vol. 3. American Mathematical
Society, 1994.

\item{[Bu1]} B. Buchberger, {\it Ein Algorithmus zum Auffinden der
Basiselemente des Restklassenringes nach einem nulldimensionalen
polynomideal}. PhD thesis, University of Innsbruck, 1965.

\item{[Bu2]} B. Buchberger, Gr\"obner bases: ~An algorithmic method
in polynomial ideal theory. In: {\it Multidimensional Systems
Theory} (Bose, N.K., ed.), Reidel Dordrecht, 1985, 184--232.

\item{[BW]} T.~Becker and V.~Weispfenning, {\it Gr\"obner Bases}.
Springer-Verlag, 1993.

\item{[CDNR]} A. Capani, G. De Dominicis, G. Niesi, and L. Robbiano,
Computing minimal finite free resolutions. {\it Journal of Pure and
Applied Algebra}, (117\& 118)(1997), 105 -- 117.

\item{[Gal]} A. Galligo, Some algorithmic questions on ideals of
differential operators. {\it Proc. EUROCAL'85}, LNCS 204, 1985,
413--421.

\item{[DGPS]} W. Decker, G.-M. Greuel, G. Pfister, and H. Sch{\"o}nemann, \newblock {\sc Singular} {3-1-3}
--- {A} computer algebra system for polynomial computations.
\newblock {http://www.singular.uni-kl.de}(2011).

\item{[Eis]} D. Eisenbud, {\it Commutative Algebra with a View
Toward to Algebraic Geometry}. GTM 150, Springer, New York, 1995. 

\item{[Fau]} J.-C. Faug\'ere. A new efficient algorithm for computing Gr\"oobner bases without
reduction to zero (F5). In: {\it proc. ISSAC'02}, ACM Press, New 
York, USA, 75-82, 2002.

\item{[Fr\"ob]} R. Fr\"oberg, {\it An Introduction to Gr\"obner Bases}. Wiley, 1997.

\item{[Kr1]} U. Kr\"ahmer,  Notes on Koszul algebras. 2010.  \\ 
 http://www.maths.gla.ac.uk/~ukraehmer/connected.pdf   
 
\item{[Kr2]} H. Kredel, {\it Solvable Polynomial Rings}. Shaker-Verlag, 1993.     

\item{[KP]} H. Kredel and M.  Pesch,  MAS, modula-2 algebra system. 1998.    
 http://krum.rz.uni-mannheim.de/mas/            
 
\item{[KR1]} M. Kreuzer, L. Robbiano, {\it Computational Commutative Algebra 1}. Springer, 2000.     
                         
\item{[KR2]} M. Kreuzer, L. Robbiano, {\it Computational Commutative Algebra 2}. Springer, 2005.

\item{[K-RW]} A.~Kandri-Rody and V.~Weispfenning, Non-commutative
Gr\"obner bases in algebras of solvable type. {\it J. Symbolic
Comput.}, 9(1990), 1--26.

\item{[Lev]} V. Levandovskyy, {\it Non-commutative Computer Algebra for
Polynomial Algebra}: {\it Gr\"obner Bases, Applications and
Implementation}. Ph.D. Thesis, TU Kaiserslautern, 2005.

\item{[Li1]} H. Li, {\it Noncommutative Gr\"obner Bases and
Filtered-Graded Transfer}. Lecture Note in Mathematics, Vol. 1795,
Springer-Verlag, 2002.

\item{[Li2]} H. Li, {\it Gr¡§obner Bases in Ring Theory}. World Scientific
Publishing Co., 2011.

\item{[Li3]} H. Li, On monoid graded local rings. {\it Journal of Pure and Applied Algebra},
216(2012), 2697 -- 2708.

\item{[Li4]} H. Li, A Note on Solvable Polynomial Algebras.
{\it Computer Science Journal of Moldova}, 1(64)(2014), 99--109. 
arXiv:1212.5988 [math.RA].

\item{[Li5]} H. Li, Computation of minimal filtered free resolutions over $\NZ$-filterd 
solvable polynomial algebras. arXiv:1401.5464 [math.RA].

\item{[LW]} H. Li and  Y. Wu, ~Filtered-graded transfer of
Gr\"obner basis computation in solvable polynomial algebras. {\it
Communications in Algebra}, 1(28)(2000), 15--32.

\item{[Sch]} F.O. Schreyer, {\it Die Berechnung von Syzygien mit
dem verallgemeinerten Weierstrasschen Divisionsatz}. Diplomarbeit,
Hamburg, 1980.

\item{[SWMZ]} Y. Sun, et al, A signature-based algorithm for computing Gr\"obner bases in
solvable polynomial algebras. In: {\it Proc. ISSAC'12}, ACM Press,
 351-358, 2012.

\end{document}